\documentclass[12pt]{article}
\usepackage[english]{babel}

\newcommand{\upto}{\nearrow}

\newcommand{\dist}{\mbox{\rm dist}}
\newcommand{\fall}{ \ \ \forall \ \ }

\newcommand{\itoinfty}{\stackrel{ i \to \infty}{\longrightarrow}}
\newcommand{\GHlim}{\stackrel{ i \to \infty}{\longrightarrow}}

\newcommand{\Riem}{{\rm Riem}}
\newcommand{\trfrriem}{    {\rm Riem \! \!\! \!\! \! {}^{^{\circ}}  \ \ }  }
\newcommand{\vol}{{\rm vol}}
\newcommand{\Happrox}{{\rm Happrox}}
\newcommand{\inj}{{\rm inj}}

\newcommand{\diam}{\mbox{\rm diam }}

\newcommand{\cig}{{\rm cig }}

\newcommand{\C}{ \mathbb C}

\newcommand{\tri}{ \triangle}
\newcommand{\curlS}{ {\cal S} }
\newcommand{\curlV}{ {\cal V} }

\newcommand{\curlM}{ {\cal M} }

\newcommand{\curlR}{ {\cal R} }

\newcommand{\curlB}{ {\cal B} }
\newcommand{\Ricci}{ {\rm Ricci}}

\newcommand{\ep}{\varepsilon}

\newcommand{\grad}{\nabla}

\newcommand{\la}{\lambda}

\newcommand{\si}{\sigma}\newcommand{\al}{\alpha}
\newcommand{\intersect}{\cap}

\newcommand{\be}{\beta}
\newcommand{\ga}{\gamma}

\newcommand{\parti}[2]{\frac{\partial #1 } {\partial #2} }
\newcommand{\lap}{\Delta}

\newcommand{\upi}{{}^{^i}\!}
\newcommand{\on}{\over}

\newcommand{\ti}{\tilde}

\newcommand{\upig}{ \upi g }
\newcommand{\upghi}{ \upi {\hat g} }

\newcommand{\uph}{^{^h}\!}

\newcommand{\upg}{^{^g}\!}
\newcommand{\upgi}{^{^{^i \!g}}\!}

\newcommand{\uptg}{^{^{\ti g}}\!}

\newcommand{\partt}{ {\partial \on {\partial t}} }

\newcommand{\Sc}{ { \rm R }}

\newcommand{\timess}{\times}

\newcommand{\br}{ {\begin{rema}}}
\newcommand{\bl}{ {\begin{lemm}} }
\newcommand{\er} { {\end{rema}}}
\newcommand{\el}{{\end{lemm}} }
\newcommand{\bc}{  {\begin{coro}} }
\newcommand{\ec}{  {\end{coro}} }

\newcommand{\gradh}{ { \uph \grad}}

\newcommand{\de}{\delta}

\newcommand{\Sp}{\mathbb{S}}
\newcommand{\N}{\mathbb{N}}
\newcommand{\R}{\mathbb{R}}
\newcommand{\disjointunion}{\dot{\cup} }
\newcommand{\GH}{ \mbox{\rm d}_{GH}  }

\newcommand{\Haus}{\mathcal{H}}

\newcommand{\T}{\mathbb{T}}

\usepackage{amsthm} 
\usepackage{amsmath}
\usepackage{amssymb} 
\usepackage{latexsym} 
\usepackage{array} 
\usepackage[T1]{fontenc} %
\usepackage{fancyhdr} 
\usepackage{xspace} %
\usepackage[dvips]{epsfig} 
\usepackage{graphicx} 

\newtheorem{defi}{Definition}[section]
\newtheorem{fact}[defi]{Fact}
\newtheorem{prob}[defi]{Problem}
\newtheorem{rema}[defi]{Remark}

\newtheorem{exam}[defi]{Example}
\newtheorem{theo}[defi]{Theorem}
\newtheorem{lemm}[defi]{Lemma}

\newtheorem{coro}[defi]{Corollary}

\newtheorem{prop}[defi]{Proposition}
\newtheorem*{thma}{Theorem A}
\newtheorem*{thmb}{Theorem B}
\newtheorem*{thmc}{Theorem C}
\newtheorem*{thmd}{Theorem D}
\newtheorem*{thme}{Theorem E}
\newtheorem*{thmf}{Theorem F}
\newtheorem*{thmg}{Theorem G}

\begin{document}

\title{ Ricci flow   of  almost non-negatively curved \\ three manifolds }
\author{Miles Simon\\
Universit{ä}t Freiburg}
\maketitle 
 \numberwithin{equation}{section}
\numberwithin{defi}{section}

\begin{abstract}

In this paper we study the evolution of 
almost non-negatively curved (possibly singular) three dimensional metric spaces by Ricci flow.
The non-negatively curved metric spaces 
which we consider arise as limits of smooth Riemannian manifolds $(M_i,\upig),$ $i \in \N,$
whose Ricci curvature
is bigger than $-{1 \on i},$ and whose diameter is less than $d_0$ (independent of $i$)
and whose volume is bigger than $v_0>0$ (independent of $i$).
We show for such spaces, that a solution to Ricci flow exists for a short time, and that the solution is smooth 
for $t>0$ and has $\Ricci \geq 0$ for $t>0.$
This allows us to classify the topological type and the differential structure of the limit manifold (in view of the theorem
of Hamilton \cite{Ha86} on closed three manifolds with non-negative Ricci curvature).

\end{abstract}

\vskip  0.07 true in
\section{Introduction and statement of results}
\vskip  0.1 true in

In Theorem \cite{Ha82} and \cite{Ha86}, Hamilton showed using the Ricci flow that

\begin{thma} \label{Ham3} {\bf (Theorem 1.2 of \cite{Ha86}) }
If $M^n,$ $n = 3 (4) $is a closed $n$-dimensional Riemannian manifold with non-negative Ricci curvature (non-negative curvature operator)
then 
$M^3$ is diffeomorphic to a quotient of $S^3$, $S^2 \timess \R,$ or $ \R^3$ by a group
of fixed point free isometries acting properly discontinuously ($M^4$ is diffeomorphic to a quotient of one
of the spaces $S^4$, $\C P^2$, $S^2 \timess S^2,$ $S^3 \timess \R^1$, $S^2 \timess \R^2$ or $\R^4$ 
by a group of fixed point free isometries acting properly discontinuously) in the standard metric.
\end{thma}

It is interesting to note that in
order to apply the theorem for $n =3$ we only require information on the Ricci curvatures (not the sectional curvatures).
The theorem implies that only certain three manifolds admit Riemannian metrics with non-negative
Ricci curvature. This is not the case for negative Ricci curvature, as proved by Lohkamp in
\cite{Lo}: he proved that every closed manifold of dimension $n \geq 3$
admits a Riemannian metric of negative Ricci curvature.

We say that a smooth family of metrics $(M,g(t))_{t \in [0,T)}$ is a solution to the Ricci flow with initial value $g_0$, or is a 
Ricci flow of $g_0$ if
\begin{eqnarray}
&&\partt g(t)= - 2 \Ricci(g(t)) \ \forall \ t \in [0,T), \cr
&&g(0) = g_0 
\end{eqnarray}
In three (and four) dimensions, there are similar results to Theorem A 
allowing some negative curvature, as the following theorems illustrate.
In Simon \cite{Si} it was shown that 

\begin{thmb}{\bf (Theorem 1.2 of \cite{Si}) }
For every $d_0,v_0 >0 $
there exists a constant $\ep = \ep(d_0,v_0) > 0 $ such that 
if 
$(M^n,g),$ $n = 2,3$ ($n =4$), is a closed Riemannian manifold with 
\begin{eqnarray}
&&\sup_{M} | \Riem | \leq 1  \label{a1} \\
&&\vol \geq v_0 > 0  \label{b1} \\
&&\diam \leq d_0 < \infty \label{c1} \\
&&\Ricci \geq -\ep \ \ ( \curlR \geq -\ep)  \label{d1}  
\end{eqnarray}
then $M$ admits a smooth metric with non-negative Ricci curvature (curvature operator:
here $\curlR $ refers to the curvature operator of $(M,g)$).
Hence, $M^n$ must be one of the possibilities listed in Theorem A.
\end{thmb}

Notice that after scaling so that $\sup | \Riem | = 1$ (which implies condition \eqref{a1}),  
\eqref{b1}  and \eqref{c1} always hold for some constants $d_0,v_0.$
Hence the condition \eqref{d1} says that $M$ is not too negatively Ricci (curvature operator) curved once we have 
scaled so that  \eqref{a1},   \eqref{b1},  and \eqref{c1} hold.

Before continuing this discussion, we would like to say a few words about the method of proof
of this theorem, as this method will be a guide 
line for the discussion and proofs which follow.

Let $d_0,v_0$ be fixed and assume that the theorem is not true (fix $n = 3$). Then there exists a sequence of Riemannian manifolds
$(M^3_i,{^i g})$ satisfying   \eqref{a1},   \eqref{b1},  and \eqref{c1}, and so that
$\Ricci \geq - {1 \on i},$ but so that $M^3_i$ {\it does not} admit a metric with $\Ricci \geq 0.$

\begin{defi}
Let $\curlB(n,d_0,-k^2_1,k^2_2,v_0)$ denote the space of smooth $n$-dimensional Riemannian manifolds, whose sectional curvatures 
are bounded from below by $-k^2_1$ bounded above by $k^2_2$, whose diameter is bounded above by $d_0$
and whose volume is bounded below by $v_0$. 
\end{defi}

It is well known that
the space $\curlB(n,d_0,-k^2_1,k^2_2,v_0)$
is precompact in the Gromov-Hausdorff space. That is, given 
a sequence of smooth $n$-dimensional Riemannian manifolds $(M^n_i,g_i)_{i \in \N} \in \curlB(n,d_0,-k^2_1,k^2_2,v_0)$
there exists a metric space $(X,d_{\infty})$ and a subsequence 
of $(M^n_i,g_i)$ (which we also call $(M^n_i,g_i)$ for ease of reading) such that
$(M^n_{i}, d(g_{i}) ) \GHlim (X,d_{\infty}),$ in the Gromov-Hausdorff sense,
where here $d(g)$ denotes the distance function (metric) $d(g):M\timess M \to \R_0^+$
arising from the Riemannian metric $g$ (see Appendix A).
The Gromov Hausdorff (space) distance between two metric spaces is defined in Appendix A.
It is a very weak measure of how close two metric spaces are to being isometric to one another.

As $\curlB(n,d_0,-k^2_1,k^2_2,v_0)$ has very well controlled geometry, it turns out that
in this case all possible limits $X$ are compact $n$-dimensional manifolds, and that $d_{\infty}$
is the distance metric arising from a $C^{1,\al}$ Riemannian metric $g_{\infty}$.
Furthermore, $X$ is diffeomorphic to $M_{i}$ for large enough $i$. 
See \cite{Pet}.
In \cite{Si} it was shown, using Ricci flow techniques, that it is possible to {\it smooth out} the
limit metric $g_{\infty}$ to obtain a smooth metric $g$ so that $(X,g)$ has non-negative Ricci curvature
(a Ricci flow of $g_{\infty}$ was constructed by taking a limit of the solutions arising from each $g_i$). 
This contradicts the assumption that $M_i$ does not admit a metric with $\Ricci \geq 0.$.
Hence, such an $\ep$ must exist.

As we mentioned above, $\curlB(n,d_0,-k^2_1,k^2_2,v_0)$ has very well controlled geometry, and so its closure does not contain 
very irregular Riemannian manifolds.

In \cite{Ya}, Yamaguchi showed that
\begin{thmc}{\bf (Corollary in the introduction of \cite{Ya})}
There exists a positive number $\ep$ such that if the curvature and diamater of a closed Riemannian three-manifold satisfy
$K_M \diam(M)^2 > -\ep$ then the possible topological type of $M$ is one of the following
\begin{itemize}
\item[{(a)}] $M$ is diffeomorphic to a torus or an infra-nilmanifold
\item[{(b)}] $M$ is up to a finite cover diffeomorphic to an $S-$bundle over $S^1$, where $S$ is one of $S^2, \R P^2, \T^2 $ 
or a Klein bottle.
\item[{(c)}] $M$ is a rational homology Sphere.
\end{itemize}
\end{thmc}

This is quite a good topological characterisation of the class of manifolds with $K_M \diam^2(M) > -\ep$ and $\ep$ small,
apart from the last possibility $(c)$, which is not very explicit. 
In a later work, Fukaya and Yamaguchi proved that for the class of Riemannian manifolds mentioned in this theorem,
the fundamental group is almost nilpotent (that is it contains a subgroup which is
nilpotent).
This allowed them to improve the classification in case (c). With this result, they obtained the (better) corollary: 

\begin{thmd}{\bf (Corollary 0.13 of \cite{FuYa})}
There exists an $\ep >0$ such that if $(M^3,g)$ is a Riemannian manifold whose 
diameter is not larger than $1,$
and has $\sec \geq - \ep,$ then a finite covering of $M$ is either
\begin{itemize}
\item homotopic  to an $S^3$ or 
\item diffeomorphic to one of
\begin{itemize}
\item[{(a)}] $T^3$
\item[{(b)}]  $S^1\timess S^2$  
\item[{(c)}] $Nil$
\end{itemize}
\end{itemize}
\end{thmd}

Hence,
using that the Poincaré Conjecture is correct (see Perelman's 
papers \cite{Pe1}, \cite{Pe2}) (that is, a homotopy $S^3$ is homeomorphic to $S^3$), we have
a good topological classification of 3-manifolds with $\sec \cdot \diam^2 \geq -\ep$ and $\ep$ small enough.

In studying the class of manifolds with $\sec \diam^2 \geq -\ep$ ($\ep$ small)
everything becomes a lot more complicated than in the situation of Theorem B above. 
There are essentially two extra problems that can now occur when studying this problem: {\it collapsing}
and {\it no upper curvature bound}. We explain this in the following.

\begin{defi}
 For $n \in \N$, $d_0 \in \R^+$ and $k \in \R$
let $\curlS(n,d_0,k)$ denote the space of smooth $n$-dimensional Riemannian manifolds of dimension $n$ 
with diameter bounded by $d_0$ and sectional curvature not less than $k.$
\end{defi}

It is well known, see for example proposition 10.7.3, and Remark 10.7.5  in \cite{BuBuIv}
that the space $\curlS(n,d_0,k)$ is precompact in the Gromov-Hausdorff space, just as $\curlB$ was.
That is, given  a sequence of smooth $n$-dimensional Riemannian manifolds 
$(M^n_i,g_i)_{i \in \N} \in \curlS(n,d_0,k),$
there exists a metric space $(X,d_{\infty})$ and a subsequence 
 of $(M^n_i,g_i)$(for ease of reading we will denote this subsequence also as $(M^n_i,g_i)$)
such that $(M_i,d(g_i))$
converges to $(X,d_{\infty})$ in the Gromov Hausdorff space.

It is however now possible that $X$ is not diffeomorphic to the $M_i$'s (for large $i$) and 
$X$ could be a manifold (with or without boundary) of lower dimension.
One may see this easily 
by considering the following example of shrinking torii

\begin{exam}
Let  
$$(M_i,g_i) = (S^1\timess \ldots\timess S^1,f_1(i)d\al^2 \oplus f_2(i)d\al^2 \ldots \oplus f_n(i) d \al^2),$$
$i \in \N$ where $d\al^2$ is the standard metric on $S^1,$ and $f_j(i) \in \R^+  \itoinfty 0,$ for all $j \in \{1, \ldots, n-m\},$ and
$ f_j(i) = 1$ for all $j \in \{ n-m +1, \ldots, n \}$
(that is, $n-m$ circles shrink to points as $i \to \infty$ and the  others stay fixed).
Then  $(M^n_i,d(g_i))$ converges in the Gromov Hausdorff space to 
$(\T^{m},d_{\T^{m}} ),$ the smooth standard $m-$dimensional torus.
\end{exam}

\begin{defi}
If $$\vol(M_i,g_i) \itoinfty 0 $$ for a sequence  of smooth Riemannian manifolds $(M_i,g_i)$
then we say that the sequence is a {\it collapsing} sequence, or that the sequence {\it collapses}.
If there exists a $v_0>0$ such that
$$\vol(M_i,g_i) \geq v_0 \ \forall \ i \in \N,$$ then we say that the 
sequence  is a {\it non-collapsing} sequence, or that the sequence {\it does not collapse}.
\end{defi}

It is also possible that the limit space $(X,d_{\infty})$ does not enjoy the regularity properties of 
the spaces occurring in the converging sequence, as one sees in the following
example.

\begin{exam}\label{coneex}
Let $(S^n_i,g_i)_{i \in \N}$ be a sequence of spheres with Riemannian metrics, where the metrics are chosen so that
\begin{itemize}
\item the sectional curvature is non-negative
\item the manifolds are becoming cone like in a fixed compact region (topologically a closed disc) as $i \to \infty$, and stay smooth away from this region.
\item the diameter is bounded above by $0 < d_0 < \infty $ and the volume bounded below by $v_0>0$ where
$d_0,v_0$ are constants independent of $i \in \N$.

\end{itemize}
Then  $(S^n_i,d(g_i))$ converges in the Gromov Hausdorff space to $(S^n,d)$, where $d$ is a (non-standard) metric on the sphere,
and there exists a Riemannian metric 
$g$ which is smooth away from the tip, induces $d$,
but cannot be extended in a $C^0$ way to the tip. It is not possible to find a $C^0$ Riemannian metric $g$ which induces $d$.
\end{exam}

So even when a sequence $(M^n_i,g_i)_{i \in \N} \in \curlS(n,d,k) $ is non-collapsed,
it is possible that the limit space $(X,d)$ may be {\it irregular} 
if the curvature of the manifolds in the sequence is not bounded from above.

In \cite{ShYa00}, further results about collapsing sequences of manifolds 
in $\curlS(n,d_0,-1)$ were proved by Shiyoa and Yamaguchi. In particular, they examined the cases that
the limit $X$ of a collapsing sequence in  $\curlS(3,d_0,-1)$ is a {\it one or two dimensional 
Alexandrov space with curvature bounded from below by $-1$} (that is, the metric space in question satisfies
certain comparison inequalities, similar to those satisfied by Riemannian manifolds 
with curvature bounded from below by $-1$: see Appendix A). 
They proved

\begin{thme}\label{cor07}{\bf (Theorems 0.2 - 0.7 of \cite{ShYa00}) }
For a given $d_0< \infty$, there exists a positive constant $\ep = \ep(d_0)>0$ satisfying the following:
If 
\begin{eqnarray}
(M,g) && \in \curlS(3,d_0,-1), \cr
\vol(M,g) && \leq \ep
\end{eqnarray}
then one of the following holds  
\begin{itemize}
\item $M$ is a graph manifold, or 
\item a finite cover of $M$ is a simply connected
manifold.
\end{itemize}
\end{thme}

Hence, using that the Poincaré Conjecture is correct  (see Perelman \cite{Pe1} , \cite {Pe2} ),  we obtain a 
good classification of three dimensional manifolds in $\curlS(3,d_0,-1)$ 
with  $\vol \leq \ep$ for $\ep = \ep(d_0)>0 $ small enough.

In \cite{ShYa05} Shiyoa and Yamaguchi continued their study of such collapsing sequences.
In their paper \cite{ShYa05} they consider the case
that the diameter of a collapsing sequence may go to infinity.
They showed

\begin{thmf}{ \bf  (Theorem 1.1 of \cite{ShYa05})}
There exist  $ 0< \ep_1 < \infty$ and $d_1 < \infty$
such that if $(M^3,g)$ satisfies $\sec  \geq -1$ and
$\vol \leq \ep_1$ then either
\begin{itemize}
\item{[case 1]}
$M$ is homeomorphic to a graph manifold or,
\item{[case 2]}
a finite cover of $M$ has finite fundamental group, and $\diam(M,g) \leq d_1$.
\end{itemize}
\end{thmf}
Hence, once again, using that the Poincaré Conjecture is correct we obtain a 
good classification of three dimensional manifolds $M^3$ with
$\sec \cdot \vol^{2 \on 3} \geq - \ep$ for $\ep$ small enough.

In this paper we consider the somewhat weaker space of {\it Riemannian manifolds with Ricci curvature bounded from below}.
\begin{defi}
For $n \in \N$, $d_0 \in \R^+$ and $k \in \R$
let $\curlM(n,d_0,k)$ denote the space of smooth $n$-dimensional Riemannian manifolds of dimension $n$ 
with diameter bounded above by $d_0$ and Ricci curvature not less than $k.$
\end{defi}

Clearly $\curlS(n,d_0,k) \subset \curlM(n,d_0,(n-1)k).$ 

In Theorem 8.6 of \cite{ChCo96}, Cheeger/Colding were able to show a similar result to that of Fukaya/Yamaguchi
concerning the Fundamental group of manifolds with almost non-negative Ricci curvature (and hence establish that
a conjecture of Gromov is correct).
They showed that
\begin{thmg}(Gromov's conjecture)
If $M$ is an  $n$-dimensional manifold with $$\Ricci \cdot \diam^2 \geq -\ep,$$ where $\ep$ (depending on $n$) is small enough,
then the fundamental group of $M$ is almost nilpotent.
\end{thmg}
The new result which allowed them to prove this theorem, is their splitting theorem, Theorem 6.64 of \cite{ChCo96},
for Gromov-Hausdorff limits of Riemannian manifolds $(M_i,g_i)$ with  $\Ricci(g_i)  \geq -\ep(i),$ where  
$\ep(i) \itoinfty 0.$

In this paper we will be chiefly concerned with metric spaces $(M^3,d_{\infty})$ which arise 
as Gromov-Hausdorff limits of non-collapsing sequences of Riemannian manifolds $(M^3_i,g_i) \in \curlM(3,d_0,-\ep(i))$ 
where $\ep(i) \to 0$ as $i \to \infty$. We show that a Ricci flow solution with initial data $(M^3,d_{\infty})$ exists for a short time (see 
Theorem \ref{appli}), and as a result, we obtain the following classification of the possible topologial types of $M$:

\begin{theo}\label{maintheo}
For all $0< v_0 < \infty , 0 < d_0 < \infty $ there exists an $\ep = \ep(v_0,d_0) >0$
such that if $(M^3,g)$ is closed and satisfies $$\vol(M,g) \geq v_0 $$ and
$(M,g) \in   \curlM(3,d_0,-\ep)$
then $M$ is diffeomorphic to a quotient of
 $S^3,$ $S^2\timess \R$ or  $\R^3$ by a group of fixed-point free isometries
acting properly discontinuously.
\end{theo}

We may write this in a more scale invariant form:

\begin{theo}
Let $d_0$ be given. There exist  $0<\ep_2 = \ep_2(d_0)< \infty$
such that if $(M^3,g)$ satisfies 
\begin{eqnarray}
&&\Ricci \cdot {\vol}^{2 \on 3}  \geq -\ep_2 \cr 
&& \diam^3 \leq  d^3_0 \cdot \vol 
\end{eqnarray}
then $M$ is diffeomorphic to a quotient of
 $S^3,$ $S^2\timess \R$ or  $\R^3$ by a group of
fixed-point free isometries acting properly discontinuously.
\end{theo}

\begin{rema}
In the case that we have a sequence of manifolds
 $(M_i,g_i)$ with $ \sec \geq - \ep(i)$, $\ep(i) \itoinfty 0,$ 
and $\vol \geq v_0,$ $\diam \leq d_0,$ then we obtain the same classification:
this is a slight improvement of the  Theorem D ( Corollary 0.13 of \cite{FuYa}), as we may remove possibility $(c)$ (that
$M$ is a nilmanfiold)
from the list of possibilities in Theorem D.
\end{rema}

As a byproduct, we note that if we combine
Theorem F and Theorem \ref{maintheo} we obtain the same classification which is implied by using Theorem F, 
without having to use that
the Poincaré Conjecture is correct.

\begin{theo}
Let $\ep_1$ and $d_1$ be as in Theorem F.
There exist  $0<\ep_3< \infty$
such that if  $(M^3,g)$ satisfies $\sec \cdot {\vol}^{2 \on 3}  \geq -\ep_3$ and 
\begin{itemize}
\item{[case 1]}
$\diam^3 \geq { 2(d_1)^3 \on \ep_1} \vol $
then $M$ is homeomorphic to a graph manifold
\item{[case 2]} $\diam^3 \leq  { 2(d_1)^3 \on \ep_1} \vol $
 then $M$  is diffeomorphic to a quotient of
 $S^3,$ $S^2\timess \R$ or  $\R^3$ by a group of fixed-point free isometries
acting properly discontinuously.
\end{itemize}
\end{theo}

\begin{proof}
Assume the theorem is false.
Then there exists a sequence of Riemannian manifolds $(M^3_i,g_i)$
with $ \sec \cdot {\vol}^{2 \on 3}  \geq -\ep_i$ for a sequence $\ep_i \in \R^+$
with $\ep_i \itoinfty 0$ and so that either  (after taking a subsequence)
\begin{itemize}
\item{case1}   $\diam^3(g_i) \geq {2(d_1)^3  \on \ep_1} \vol(g_i) $   but  $M_i$ is not homeomorphic
to a graph manifold, or
\item{case2}  $\diam^3(g_i) \leq {2(d_1)^3  \on \ep_1} \vol(g_i) $ but $M_i$ is not one of the possibilities
listed in case 2
of the statement of this theorem.
\end{itemize}
In case 1, we scale all $(M_i,g_i)$ so that $\vol(g_i) = \ep_1.$
Then $\diam(g_i) \geq  2 d_1,$ and $\sec \geq -1$ for $i$ sufficiently large, 
and so we may apply Theorem F to obtain a contradiction.

In case 2, let us scale all $(M_i,g_i)$ so that $\vol(g_i) = 1.$ We then have
$\diam(g_i) \leq  { 2 d_1 \on (\ep_1)^{1 \on 3}},$ and 
 $\sec  \geq  - \ep_i$
and so we may then apply the Theorem \ref{maintheo} to obtain a contradiction.
\end{proof}

In this paper we will cheifly be
concerned with metric spaces $(M^3,d_{\infty})$ which arise 
as Gromov-Hausdorff limits of non-collapsing sequences of Riemannian manifolds $(M^3_i,g_i)  \in \curlM_i(3,d_0,-\ep(i))$ 
where $\ep(i) \to 0$ as $i \to \infty.$
In particular, we wish to flow such metric spaces  $(M^3,d_{\infty})$ by Ricci Flow.
As we saw in the previous section (see Example \ref{coneex}) such limits can be quite irregular (not even $C^0$).
Nevertheless, they will be Alexandrov Spaces of curvature bounded from below by $0$, and so do carry
some structure (see Appendix A).
In order to flow  $(M^3,d_{\infty})$ we will flow each of the $(M^3_i,g_i)$ and then take a Hamilton limit of the
solutions (see \cite{Ha95a}).
The two main obstacles to this procedure are:
\begin{itemize}
\item It is possible that the solutions $(M_i,g_i(t))$  are defined only for $t \in [0,T_i)$ where $T_i \to 0$ as $i \to
\infty$.
\item In order to take this limit, we require that each of the solutions satisfy uniform bounds of the form
$$\sup_{M_i} |\Riem(g_i(t))| \leq c(t), \ \forall \ t \in (0,T) ,$$
for some well defined common time interval 
$(0,T)$ ($c(t) \to \infty$ as $t \to 0$ would not be a problem here). Furthermore they should all
satisfy a uniform lower bound on the injectivity radius of the form
$$\inj(M,g_i(t_0)) \geq \si_0 > 0 $$ for some $t_0 \in (0,T).$
\end{itemize}
As a first step to solving these two problems,  
in Lemma \ref{globaltimelemma} of Section \ref{chapbound} we see that a (three dimensional) smooth 
solution to the Ricci flow $(M,g(t))_{t \in [0,T)}$ cannot become singular at time $T$ as
long as $\Ricci \geq -1,$ the diameter remains bounded (by say $d_0$)
and the volume stays bounded away from zero (say it is bigger than $v_0$).
Furthermore, a bound of the form
$$| \Riem(g(t))|  \leq {c_0(d_0,v_0) \on t} \forall t \in [0,T) \cap [0,1]$$
for such solutions is proved: that is, the curvature of such solutions is quickly smoothed out.

In  Lemma \ref{globaltimecor} we present an application of the proof of \ref{globaltimelemma}.
Notice that Proposition
11.4 of \cite{Pe1} for the three dimensional case implies the Lemma
\ref{globaltimecor}. Perelman's method of proof is somewhat different from that used in  Lemma \ref{globaltimecor}. 

Section 4 is concerned with proving (for an arbitrary three dimensional solution to the Ricci flow)
lower bounds for the Ricci curvature of the evolving metric, which depend on 
\begin{itemize}
\item the bound from below for the Ricci curvature of the initial metric
\item the scalar curvature of the evolving metric.
\end{itemize}
One of the major applications is: 
if $(M,g_0)$ satisfies $\Ricci(g_0) \geq -\ep_0$ ($\ep_0$ small enough) and the solution satisfies
$\Sc(g(t)) \leq {c_0 \on t},$ then
$$ \Ricci(g(t)) \geq - 2 c_0\ep_0 \forall t \in (0,T_*) \cap (0,T) $$
for some universal constant $T_* = T_* >0$ ( $(0,T)$ is the time interval for which the solution
is defined). See Lemma \ref{seclemma2}.

In Section 5, we consider smooth solutions to the Ricci flow which satisfy
\begin{eqnarray}
&&\Ricci(g(t)) \geq - c_0  \label{intro1} \\
&& |\Riem(g(t))| t \leq c_0 \label{intro2}  \\
&& \diam (M,g_0) \leq d_0 \label{intro3}
\end{eqnarray}
In Lemma \ref{diamlemm}, well known bounds on the evolving distance for a solution to the Ricci flow are proved
for such solutions.

We combine this Lemma with some results on Gromov-Hausdroff convergence to show (Corollary \ref{volcorr})
that such solutions can only lose volume at a controlled rate.

In Section 6 we show (using the a priori estimates from the previous sections) that
a solution to the Ricci flow of $(M,d_\infty)$ exists, where  $(M,d_\infty)$ is the Gromov-Hausdorff limit
as $i \to \infty$ of $(M_i, d(g_i))$ where the 
$(M_i, g_i)$ satisfy 
\begin{eqnarray*}
&&\Ricci(g_i) \geq - \ep(i)  \cr
&& \vol(M,g_i) \geq v_0 \cr
&& \diam (M,g_0) \leq d_0. 
\end{eqnarray*}

More explicitly we prove:
\begin{theo}
Let $(M_i, \upi g_0)$ be a sequence of closed three (or two) manifolds satisfying
\begin{eqnarray*}
&&\diam(M_i,\upi g_0) \leq d_0 \cr
&&\Ricci(\upi g_0) (\sec(\upi g_0) )  \geq  - \ep(i) {\upi g}_0 \cr
&&\vol(M,\upi g_0) \geq v_0 >0,
\end{eqnarray*}
where $\ep(i) \to 0,$ as
$i \to \infty.$ 
Then there exists an $S = S(v_0,d_0) >0$ and $K = K(v_0,d_0)$ such that the
maximal solutions $(M_i,\upig(t))_{t \in [0,T_i)}$ to Ricci-flow satisfy
$T_i \geq S,$
and $$\sup_{M} |\Riem(\upig(t))| \leq {K \on t},$$ for all
$t \in (0,S)$.
In particular the Hamilton limit solution 
$(M,g(t))_{t \in (0,S)} = \lim_{i \to \infty} (M_i,\upig(t))_{t \in (0,S)}$ (see \cite{Ha95a})
exists and satisfies
\begin{eqnarray}
&& \sup_{M} |\Riem(g(t))| \leq {K \on t} \cr
&& \Ricci(g(t)) \geq 0  \ \ (\sec(g(t)) \geq 0),
\end{eqnarray}
for all $t \in (0,S)$ and $(M,g(t))$ is closed. Furthermore 
\begin{eqnarray}
\GH( (M, d(g(t))  ) , (M, d_\infty) ) \to 0 
\end{eqnarray}
as $t \to 0$ where
$(M,d_\infty) = \lim_{i \to \infty} (M_i,d(\upi g_0) )$
(the Gromov-Hausdorff limit). Hence, if $M = M^3$, then $M^3$ is diffeomorphic to
a quotient of one of $S^3$,$ S^2 \timess \R $ or $ \R^3$ by group
of fixed point free isometries acting properly discontinuously. 
\end{theo}

The theorem which is essential in constructing such a solution is:
\begin{theo}
Let $M$ be a closed three (or two) manifold satisfying
\begin{eqnarray}
&&\diam(M,g_0) \leq d_0 \cr
&&\Ricci(g_0) \ \ \ (\sec(g_0) )  \geq  - \ep g_0 \cr
&&\vol(M,g_0) \geq v_0  >0,
\end{eqnarray}
where $\ep \leq  {1 \on 10 c^2}$ and $c = c(v_0,d_0) \geq 1$ is the constant from Lemma 
\ref{globaltimelemma}
Then there exists an $S = S(d_0,v_0) >0$ and $K = K(d_0,v_0)$ such that the
maximal solution $(M,g(t))_{t \in [0,T)}$ to Ricci-flow satisfies
$T\geq S,$
and $$\sup_{M} |\Riem(g(t))| \leq {K \on t},$$ for all
$t \in (0,S)$.
\end{theo}

Appendix A contains definitions, results and facts about Gromov-Hausdorff space, which
we require in this paper. 

In Appendix B we define $C-$ essential points, and $\de$-like necks, and consider
discuss $0$-like necks in the three dimensional case.

A proof of the (well known) Lemma \ref{diamlemm} is contained in Appendix C.

Appendix D is a description of the Notation used in this paper.

\section{ Bounding the blow up time from below using bounds on the geometry.}\label{chapbound}
\setcounter{equation}{0}
\setcounter{defi}{0}
 \numberwithin{defi}{section} 
\numberwithin{equation}{section}
%


An important property of the Ricci flow is that:\hfill\break
if certain geometrical quantities are controlled (bounded) on a half open finite time interval $[0,T)$, then
the solution does not become singular as $t \upto T$ and may be extended to a solution defined on the time interval
$[0,T + \ep)$ for some $\ep >0.$
We are interested in the question:

\begin{prob}
What elements of the geometry need to be controlled, in order to guarantee that a solution does not become singular?
\end{prob}

In \cite{Ha82}, it was shown that
for $(M,g_0)$ a closed smooth Riemannian manifold, the Ricci flow equation
\begin{eqnarray}
&&\partt g = - 2 \Ricci(g) \label{Ricci}\cr
&&g(\cdot,0) = g_0,
\end{eqnarray}
always has a solution
$(M,g(t))_{t \in [0,T)}$
for a short time. It was also shown that two such solutions  defined on the same time interval must agree, if there
initial values agree.
Furthermore, for each smooth, closed $(M,g_0)$ there exists a maximal time interval $[0,T_{Max})$ ($T_{Max} >0$) for which, 
there exists a solution $(M,g(t))_{t \in [0,T_{Max})}$ to \eqref{Ricci}, and if $T_{Max} < \infty$ then there
is no solution  $(M,g(t))_{t \in [0,T_{Max}+ \ep)}$ to \eqref{Ricci} (for any $\ep >0$ ). Such a solution
 $(M,g(t))_{t \in [0,T_{Max}) }$ is called a {\it maximal solution}.
\begin{defi}(Maximal Solutions)
Let $(M,g(t))_{t \in [0,T)}$ be a solution to Ricci flow. We say that the solution {\it blows up} at time $T$
if 
\begin{eqnarray}
\sup_{M\timess [0,T)} |\Riem| = \infty,
\end{eqnarray}
\end{defi}

It was also shown in \cite{Ha82} that
\begin{lemm}
Let $(M,g(t)_{t \in [0,T)}$ be a closed, smooth solution to Ricci flow, with $g(0) = g_0$
and $T < \infty,$ with
\begin{eqnarray}
\sup_{M\timess [0,T)} |\Riem| < \infty.
\end{eqnarray}
Then, for some $\ep >0,$ there exists a solution $(M,g(t))_{t \in [0,T + \ep)}$, with $g(0) = g_0$.
\end{lemm}

So we see that a bound on the supremum of the Riemannian curvature (that is, {\it control}
of this geometrical quantity) on a finite time interval $[0,T)$
guarantees that this solution does not become singular as $t \upto T$.
In the following lemma, we present 
other bounds on geometrical quantities which guarantee that a solution to the Ricci flow does
not become singular as $t \upto T$.

\begin{lemm}\label{globaltimelemma}

Let $(M^3,g(t))_{t \in [0,T)}$, $T \leq 1$ be an arbitrary smooth solution to Ricci flow  ($M^3$ closed ) satisfying
\begin{eqnarray}
\Ricci(g) & \geq & -1, \cr
\vol(M,g) & \geq & v_0 > 0 \cr
\diam(g) & \leq & d_0 < \infty \label{volumdiam}
\end{eqnarray}
 for all $t \in [0,T).$
Then there exists a $c = c(d_0,v_0),$
such that
\begin{eqnarray*}
\Sc(g(t)) t \leq c
\end{eqnarray*}
for all $t \in [0,T).$
In particular,  $(M^3,g(t))_{t \in [0,T)}$ is not maximal.
\end{lemm}

\begin{coro}\label{globaltimecor}
Let $(M^3,g(t))_{t \in [0,T)}$ be an arbitrary smooth solution to Ricci flow  satisfying
\begin{eqnarray}
\Ricci(g) & \geq & -1, \cr
\vol(M,g) & \geq & v_0 > 0 \cr
\diam(g) & \leq & d_0 < \infty \label{volumdiam2}
\end{eqnarray}
 for all $t \in [0,T).$
Then there exists a $c = c(d_0,v),$
such that
\begin{eqnarray*}
\Sc(g(t))  \leq c \max( {1 \on t}, 1 )
\end{eqnarray*}
for all $t \in [0,T).$
In particular,  $(M^3,g(t))_{t \in [0,T)}$ is not maximal.
\end{coro}

The proof of the corollary is a trivial iteration argument.

\begin{proof} (of the corollary)
Fix $t_0 \in [0,T).$
We wish to show that $$\Sc(g(t_0))  \leq c \max( {1 \on t_0}, 1 ).$$
If $t_0 \leq {1 \on 2}$ then we apply Lemma \ref{globaltimelemma}.
If ${(N + 1) \on 2} > t_0 \geq {N \on 2},$ ($N \in \N$) then we apply Lemma \ref{globaltimelemma}
to the solution $(M,g( {(N-1) \on 2 } + t) )_{t \in [{1 \on 2},1) } $ of Ricci flow
(notice that  ${(N-1) \on 2 } + t = t_0$ implies that 
$1 > t \geq {1 \on 2}$).
\end{proof}

We now prove Lemma \ref{globaltimelemma}.

\begin{proof}

Assume to the contrary that there exist solutions $(M_i, {\upig(t)})_{t \in [0,T_i)}$, $T_i \leq 1$ to Ricci flow
such that
 \begin{eqnarray}
 \sup_{(x,t) \in M_i\timess (0,T_i)} {\upi \Sc  }(x,t)t  \itoinfty \infty, 
\end{eqnarray} 
or there exists some $j \in \N$ with
\begin{eqnarray}
 \sup_{(x,t) \in M_j\timess (0,T_j)} {{}^j \Sc  }(x,t)t  =  \infty, 
\end{eqnarray} 
where ${\upi \Sc  }:= \Sc({\upig}).$
It is then possible to choose points $(p_i,t_i) \in M_i\timess [0,T_i)$ (or in  $M_j \timess [0,T_j)$: in this case we
redefine $M_i =M_j$ and $T_i = T_j$ for all $i \in \N$ and hence we do not need to treat this case separately ) such that
\begin{eqnarray}
\Sc(p_i,t_i)t_i  = \sup_{(x,t) \in M_i\timess (0,t_i]} {\upi \Sc  }(x,t)t  \itoinfty \infty.
\end{eqnarray}

Define
\begin{eqnarray}
  { \upghi}(\cdot,\hat t) := c_i {{\upig}}(\cdot,t_i +{\hat t \on c_i}),
\end{eqnarray}
 where
$c_i := {\upi \Sc  }(p_i,t_i).$
This solution to Ricci flow is defined for $ 0 \leq t_i +{\hat t \on c_i}< T_i,$
that is, at least for $0 \geq \hat t > -t_i c_i.$
Let $A_i := t_i c_i.$ Then the solution  ${ \upghi}(\hat t)$ is defined at least for
$  \hat t \in (-A_i, 0 ).$
By the choice of $(p_i,t_i)$ we see that the solution is defined
for $\hat t > -A_i =  -t_i c_i =  -t_i{\upi \Sc }(p_i,t_i) \itoinfty - \infty.$
Since $t_i \leq T_i  \leq 1,$ we also have 
\begin{eqnarray} 
c_i \itoinfty \infty , \label{cieqn}
\end{eqnarray}
in view of the fact that
 $$ t_i c_i= t_i{\upi \Sc }(p_i,t_i) \itoinfty \infty .$$

Furthermore, letting $s(\hat t,i) := t_i +{\hat t \on c_i},$
where $  -A_i <  \hat t \leq 0$ we have
\begin{eqnarray}
{\upi \hat \Sc}(\cdot, \hat t) &=& {1 \on c_i} {\upi{\Sc}}( \cdot,s(\hat t ,i)) \\
&=& { \upi {  \Sc}(\cdot,s) \on \upi {  \Sc}(p_i,t_i) } \cr
&=&  { \upi{ \Sc}(\cdot,s)s \on \upi{  \Sc}(p_i,t_i) t_i } {t_i \on s} \cr
&\leq&   {t_i \on s}  \cr
&=&  {t_i \on  t_i + {\hat t \on c_i} } \itoinfty 1.  \label{fidely}
\end{eqnarray}
in view of the definition of $(p_i,t_i)$, and $0 \leq s \leq t_i$
(follows from the definition of $s$ and the fact that $\hat t \leq 0$), and
\eqref{cieqn}.
Due to the conditions \eqref{volumdiam} 
we see that  there exist $l= l(v_0,d,n),$ and $ \ep =\ep(v_0,d,n),$ such that
 \begin{eqnarray}
l \geq {\vol(B_r(p),\upig(t)) \on r^3} \geq \ep \ \forall \ r \leq \diam(M_i,{\upig}(t)),
\end{eqnarray} 
(in view of the Bishop Gromov comparison principle)
which implies the same result for any rescaling of the manifolds.
Notice that the conditions \eqref{volumdiam} imply that 
\begin{eqnarray}
\diam(M,\upig(t)) \geq d_1(n,v_0) >0\label{diaminf}
\end{eqnarray}
for some  $\infty > d_1(n,v_0) > 0.$
 Otherwise, assume $\diam(M,\upig(t)) \leq d_1$ for some small $d_1$, then
$\vol(M,\upig(t)) \leq c(n) d_1^3 \omega_n $ (Bishop-Gromov comparison principle), and hence
$ \vol(M,\upig(t))  <   v_0 $ if  $d_1$ is too small, which would be a contradiction.
Hence,
$ \diam(M,  \upghi(0)) \itoinfty \infty,$
in view of the inequalities \eqref{diaminf} and \eqref{cieqn}.
Now using 
\begin{eqnarray}
l \geq {\vol(B_r(p),\upghi(t) ) \on r^3} \geq \ep_0, \forall r \leq \diam(M_i,\upghi(\hat t)), \label{volinit}
\end{eqnarray}
we obtain a bound on the injectivity radius from below, in view
of the theorem of Cheeger/Gromov/Taylor, \cite{CGT}
(the theorem of Cheeger/Gromov/Taylor says that for a complete Riemannian manifold  $(M,g)$ with
$| Riem| \leq 1$, we have
$$ \inj(x,g) \geq r { \vol(g, B_r(x)) \on   \vol(g, B_r(x)) + \omega_n\exp^{n-1} },$$
for all $r \leq {\pi \on 4}.$
In particular, using that $\diam(M,g) \geq d_1 >0$ and $|\Riem| \leq c$ (see [i] below) 
for the Riemannian manifolds in question, we obtain
\begin{eqnarray}
 \inj(x,g) \geq \ep {s^{n+1} \on l s^n +  \omega_n\exp^{n-1}} \geq
  c^2(d_0,v_0,n) > 0 
\end{eqnarray}
for $s = \min( (\omega_n \exp^{n-1})^{1 \on n},\diam(M,g),{\pi \on 4} )$).

This allows us to take a pointed {\it Hamilton limit} (see \cite{Ha95a}), which leads to a Ricci flow solution
$(\Omega,o, g(t)_{t \in (- \infty,\omega) }),$
with $R \leq R(o,0) = 1,$ and $\Ricci \geq 0,$ $\omega > 0$
(at $t = 0,$ as explained below, the full Riemannian curvature tensor of $\upghi(0)$is bounded by $c(3)$ 
and so clearly each solution lives at least to a time $\omega >0$ independent of $i$ ).
More Precisely: 
\begin{itemize}
\item{[i]} The bound from below on the Ricci curvature, and the bound from above on the scalar curvature imply
that the Ricci curvatures are bounded absolutely by the constant $5$ for $i$ big enough.
In three dimensions, bounds from above and below on the Ricci curvatures imply bounds from above and below on
the sectional curvatures and hence on the norm of the full Riemannian curvature tensor. This, together
with the bound from below on the injectivity radius, allows us to a take a Hamilton limit of these Ricci flows.
\item{[ii]} In fact the limit solution satisfies $\sec \geq 0,$ which can be seen as follows:
Each rescaled solution $\upghi$ is defined on $M_i \timess [-A_i,\omega]$ where $A_i \itoinfty \infty.$
They also each satisfy $\sec \geq -2$ and $|\Riem| \leq c(n)$
for all $t \in (-S,0)$ for any fixed $S$ and  all $i$ big enough, in view of \eqref{fidely} and
$\Ricci \geq -1.$ 

Let us translate in time by $S$ , so that these solutions are defined on 
 $M_i \timess [-A_i + S ,S]$ and satisfy  $\sec \geq -2$ and $|\Riem| \leq c(n)$ on $(0,S)$ 
(for $i$ big enough).Without loss of generality, we assume that  $\sec \geq -1$.
We then use the improved pinching result of Hamilton \cite{Ha99} (see also \cite{Iv}):
\begin{theo}
Let $g(t)$ be a solution to Ricci flow defined on $M \timess [0,T)$, $M$ closed. 
Assume at $t = 0$ that the eigenvalues 
$\al \geq \be \geq \ga$ of the curvature operator at each point are bounded below by
$ \ga \geq -1$. The scalar curvature is their sum $\Sc = \al + \be + \ga,$ and $X := -\ga$.
Then at all points and all times we have the pinching estimate
$$ \Sc \geq X [ \log X + \log(1 + t) -3 ],$$
whenever $X >0.$
\end{theo}

Notice that this estimate is also valid for the translated limit solution (defined on $[0,S)$) as 
it is valid for each $i$ and the scalar curvature and $X$  converge as $i \to \infty$ 
to the corresponding quantities of the translated (by $S$) limit solution.

Let $\de >0$ be any arbitrary small constant. Assume there exists
$(x,t) \in \Omega \timess [{S \on 2},S)$ such that
$X(x,t) \geq \de$.
Then we get
\begin{eqnarray}
\log (\de) \leq \log X(x,t)  
& \leq & {\Sc(x,t) \on \de } -\log(1+t) + 3  \cr
   & \leq & {c(n) \on \de} -\log(1+ {S \on 2}) + 3 
\end{eqnarray}
which is a contradiction for $S$ big enough.
Hence our initial limit solution (without any translations in time) has
$X(x,0) \leq \de.$
As $\de$ was arbitrary we get $X(\cdot,0) \leq 0.$
A similar argument shows $X \leq 0$ everywhere.
That is, the limit space satisfies $\sec \geq 0 \forall t \in (-\infty,0)$
\end{itemize}

The volume ratio estimates 
\begin{eqnarray}
l \geq {\vol(B_r(p)) \on r^3} \geq \ep_0 \forall r >0, \label{vol}
\end{eqnarray}
 are also valid for $(\Omega,g),$
as these estimates are scale invariant, and $\diam(\Omega,g) = \infty.$ 
At this point we could apply Proposition 11.4 of \cite{Pe1} to obtain a contradiciton.
We prefer however to introduce an alternative method to Perelman in order to obtain a contradiction
(this method may be of independent interest).
We now consider the following two cases.  
\begin{itemize}
\item[(case 1)]  $\sup_{\Omega\timess ( -\infty, 0]} |t| \Sc = \infty, $
\item[(case 2)]  $ \sup_{\Omega\timess ( -\infty, 0]}  |t| \Sc < \infty.$
\end{itemize}

\begin{itemize}
\item[(case 1)] 
In the first case, in view of \cite{ChKn}, chapter 8, section 6, we may assume w.l.o.g.
that there exists a solution $(\Omega,o, g(t)_{t \in (- \infty,\infty)}),$
with 
\begin{eqnarray}
\sup_{\Omega\timess ( -\infty, \infty)} |\Sc(t)| \leq 1 = |\Sc(0,o)| . \label{grsol} 
\end{eqnarray}
Note: we must slightly modify the argument there, by replacing $\Riem$ with $\Sc$ wherever it appears.
We also use the fact (as mentioned above) that
$|\Riem| \leq c(3) \Sc$ in the case that $\Ricci \geq 0$ (in dimension three) and that our scale invariant
volume estimate  \eqref{vol} remains true for any rescalings of our solution: these two facts ensure that 
in the rescaling process of the argument in \cite{ChKn}, chapter 8, section 6,
an injectivity radius estimate is satisfied, and that the limit solution is well defined.

\item[(case 1.1)] the sectional curvature is everywhere positive.
\item[(case 1.2)] there exists $(p_0,t_0) \in \Omega\timess ( -\infty, \infty),$
and $v_{p_0}, w_{p_0} \in T_{p_0} \Omega $ with
$\sec(p_0,t_0)(v_{p_0},w_{p_0}) = 0.$

First we consider case 1.1 
\item[(case 1.1)] 
This means $\Omega$ is diffeomorphic to  $\R^3$ in view of the soul theorem (see \cite{ChEb}, Chapter 8 ) and in particular,
$\Omega$ is simply connected.
We may then apply the gradient soliton theorem of Hamilton \cite{Ha93} which implies, in view of
\eqref{grsol}, that $(\Omega,g(t))_{t \in  ( -\infty, \infty)}$ is a gradient soliton.
We may then, using the dimension reduction theorem of Hamilton, Theorem  22.3 of \cite{Ha95}, take a Hamilton limit of 
rescalings of this solution, to obtain
a new solution,   $(\R\timess N ,dx^2 \oplus\ga(t))_{t \in  ( -\infty, \infty)},$ 
or a quotient thereof by a group of fixed-point free isometries acting properly discontinuously,
where $dx^2$ is the standard metric on $\R$ , and  $(N, \ga(t))_{t \in  ( -\infty, \infty)}$
is a solution to the Ricci flow, $N$ is a surface, and
$\Sc(\cdot,t) >0,$ on $N$. In the case that we have a 
quotient of  $(\R\timess N ,dx^2 \oplus\ga(t))$ then we notice that 
$(\R\timess N ,dx^2 \oplus\ga(t))$ still satisfies 
\eqref{vol}
(the bound from below follows as the Riemannian covering map $f:(\R\timess N ,dx^2 \oplus\ga(t)) \to (\Omega,g(t)) $ 
is a Riemannian-Submersion, 
and the bound from above follows in view of the Bishop-Gromov comparison principle)
and so, without loss of generality, we may assume that we do not have a quotient.
If $N$ is compact, then
$( \R\timess N, dx^2 \oplus \ga),$ 
does not satisfy the estimates \eqref{vol}, 
and so we obtain a contradiction. So w.l.o.g. we may assume that $N$ is non-compact.
Now we break this up into two cases:
\item[(case 1.1.1)]
 $\sup_{N\timess ( -\infty, \infty)} |t| |\Sc(t)| = \infty, $
and 
\item[(case 1.1.2)]
 $\sup_{N\timess ( -\infty, \infty)} |t| |\Sc(t)| < \infty. $

First we handle
\item[(case 1.1.1)]
Once again, w.l.o.g (\cite{ChKn} chap.8, sec.6), we may assume
 $$\sup_{N\timess ( -\infty, \infty)}  \Sc \leq 1 = \Sc(o,0). $$
$\Sc(t) >0,$ and $N$ non-compact  implies $N$ is diffeomorphic to $\R^2,$ which is simply connected.
We may then use the gradient soliton theorem of Hamilton, \cite{Ha93} , to obtain that
$(N, \ga),$ is a gradient soliton, which implies (thm. 26.3, \cite{Ha95} ), that
$(N,\ga)$ is the cigar $(\Sigma,\cig)$.
But $( \R\timess \Sigma, dx^2 \oplus \cig)$ do not satisfy the estimates \eqref{vol}, 
and so we obtain a contradiction.
\item[(case 1.1.2)]$\sup_{N\timess ( -\infty, \infty)} |t| |\Sc(t)| < \infty. $
Hamilton, Thm. 26.1 of \cite{Ha95} implies that
$(N,\ga) = (\Sp^2 \mbox{ \rm or } \R^2, \ga)/ \Gamma,$
where $\ga$ is the standard solution on $S^2$ or $\R^2,$ and
$\Gamma$ is a finite group of isometries acting without fixed points on the standard
$S^2$ or standard $\R^2.$  $(\R^2,\ga)$ cannot occur, since the surface should satisfy $\Sc(t) > 0$ everywhere (the standard 
$(\R^2,\ga)$ is flat).
But then $N$ is compact, and
$( \R\timess N, dx^2 \oplus \ga),$ does not satisfy the estimates \eqref{vol}, 
and once again we obtain a contradiction.
\item[(case 1.2)] there exists $(p_0,t_0) \in \Omega\timess ( -\infty, \infty),$
and $v_{p_0}, w_{p_0} \in T_{p_0} \Omega $ with
$\sec(p_0,t_0)(v_{p_0},w_{p_0}) = 0.$ Then the maximum principle applied to the evolution equation of the curvature operator,
implies that  $(\Omega,o, g(t))_{t \in (- \infty,\infty)} = (\R\timess N, dx^2 \oplus \ga(t))_{t \in (- \infty,\infty)},$ 
or a quotient thereof by a group of isometries(see \cite{Ha86}, chapter 9)
and $\sup_{N\timess (-\infty, \infty)} \Sc(t) \leq 1 = \Sc(o,0).$
Without loss of generality, we may assume that we don't have a quotient, as explained
in case 1.1.
$\Sc(t) >0,$ implies $N$ is diffeomorphic to $\Sp^2 / \Gamma $ or $\R^2.$
In the case that $N$ is diffeomorphic to $\Sp^2 / \Gamma$, we obtain a contradiction, as then
$(\Omega,g)$ does not satisfy \eqref{vol}.
So w.l.o.g. $N$ is diffeomorphic to $\R^2,$ in particular $N$ is simply connected.
We may use the gradient soliton theorem of Hamilton \cite{Ha93}, to get
that $(N,\ga)$ is a soliton and it must be the Cigar, in view of theorem
26.3 of Hamilton \cite{Ha95}.
This leads to a contradiction as then $(\Omega,g)$ does not satisfy \eqref{vol} (similarly for the covering case).

\item[(Case 2)]  $ B:= \sup_{\Omega\timess ( -\infty, 0]}  |t| |\Riem(t)| < \infty.$
\item[(Case 2.1)]  The asymptotic scalar curvature ratio $A = \limsup_{s \to \infty} \Sc s^2  = \infty.$
Then we use the dimension-reduction argument of Hamilton (see Lemma 22.2 of \cite{Ha95} and the argument directly after 
the proof of Lemma 22.2) to obtain a new solution
 $ (N\timess \R , \ga \oplus dx^2)$ or a quotient thereof by a group of isometries 
where $(N,\ga)$ is a solution to Ricci flow defined on  $( -\infty, T]$ ($T >0$) (note, our injectivity radius estimate is still
valid in view of the volume ratio estimate  \eqref{vol} which survives into every limit).
If $N$ is compact then we obtain a contradiction to \eqref{vol}.
So we may assume that $N$ is non-compact.
We then consider the cases
 $\sup_{N\timess ( -\infty, \infty)} |t| |\Sc(t)| = \infty, $
and 
 $\sup_{N\timess ( -\infty, \infty)} |t| |\Sc(t)| < \infty. $
Then, using the exact same arguments as in Cases 1.1.1 and 1.1.2, we obtain a contradiction.   
 \item[(Case 2.2)] The asymptotic scalar curvature ratio $A = \limsup_{s \to \infty} \Sc s^2  < \infty.$
Remember that the  asymptotic scalar curvature ratio is a constant in time for ancient solutions which have
bounded curvature at each  time and non-negative curvature operator. $A$ is also independent of which origin we choose:
see theorem 19.1 \cite{Ha95}.

Now we use another splitting argument of Hamilton (see Theorem 24.7 of \cite{Ha95} for the compact version of this argument).
\item[(Case 2.2.1)]
There exists a $C >0$, s.t., for all $\tau \in (- \infty,0),$ for all $\de \in (0,1),$
there exists $(x,t)\in {\Omega\timess ( -\infty, \tau)}$ such that $(x,t)$ is a $C$-essential $\de$-necklike point
(see Appendix  B).
Let $\{\de_i\}_{i \in \N}$ be a positive sequence,$\de_i \itoinfty 0,$ and
let $(x_i,t_i)$ be chosen so that $(x_i,t_i)$ is an $C$-essential $\de_i$-necklike point.
Assume $\theta_i$ is a unit 2-form on $T_{x_i} \Omega$ with
$$ |\Riem(x_i,t_i) - \Sc(x_i,t_i)(\theta_i \otimes \theta_i)| \leq  \de_i |\Riem|(x_i,t_i).$$
Let $ \upig(x,t) = {1 \on |t_i| }g(x,t_i + t |t_i| ).$
Then 
\begin{eqnarray}
 \upgi |\upi \Riem ( x, t)| &=& |t_i| \upg |\Riem ( x, t_i + t|t_i|)| \cr
  & = &  |t_i| {\upg |}\Riem ( x, (t-1)|t_i|)|\cr
&=& {|(t-1)|t_i|| \on | 1-t| } {\upg |}\Riem ( x, (t-1)|t_i|)| \cr
&\leq& { B \on  | 1-t| } \leq 2 B 
\end{eqnarray}
for $t \leq {1\on 2}.$
Notice that 
\begin{eqnarray}t_i + {1 \on 2}|t_i| = t_i  - {1 \on 2} t_i =  {1 \on 2} t_i< 0
\end{eqnarray} 
and so $\upig(t)$ is defined for (at least) $-\infty < t \leq {1 \on 2}$. 
Furthermore, 
\begin{eqnarray}
\upgi |\upi \Riem ( x_i, 0)| = |t_i| \upg |\Riem ( x_i, t_i)| \geq C > 0, \label{C_i}
\end{eqnarray} 
since $(x_i,t_i)$ is $C$-essential.
Set $$\psi_i := {1 \on |t_i|} \theta_i.$$ 
$\psi_i $ is then a unit two form on $T_{x_i} \Omega$ with respect to $g^i(x,0).$
Then
$$ {\upgi |}\upi \Riem(x_i,0) - \upi \Sc(x_i,0)(\psi_i \otimes \psi_i)| \leq  \de_i B.$$
Now taking a Hamilton pointed limit (our injectivity radius estimate is still valid) we obtain a solution
$(\ti {\Omega},\ti g),$ defined for $t \leq {1 \on 2}$ with 
$$ {\uptg |} \ti{\Riem}(o,0) - {\ti \Sc}(o,0)(\psi \otimes \psi )| \leq  0.$$
where $\psi$ is a unit two form defined on $T_o M,$
$\psi = \lim_{i \to \infty} \psi_i.$
More precisely this $\psi$ is obtained (in coordinates) as $\psi^{\al \be}(o):= \lim_{i \to \infty} \parti{F_i^{\al} }{x^k}(o)
\parti{F_i^{\be} }{x^l}(o) \psi_i^{kl}(x_i),$ where
$F_i: B_i(o) \to U_i \subset M_i,$ $F_i(o) = x_i$ are the diffeomorphisms occurring in the Hamilton limit process.
 
Furthermore $\Sc(o,0) \geq C >0,$ (in view of \eqref{C_i}) 
which implies (in view of the strong maximum principle applied to 
the evolution equation for $\Sc$) that $\Sc >0$.
Hence, due to the maximum principle,
$(\ti {\Omega},\ti g) = (N\timess \R, \ga \oplus dx^2),$ or a quotient thereof by a group of isometries,
where $(N,\ga)$ is a solution to the Ricci flow 
(see Appendix B for a more detailed
explanation of this fact).
If $N$ is compact we obtain a contradiction to the 
volume ratio estimates.
If $N$ is non-compact, then we argue exactly as in cases 1.1.1 and 1.1.2 to obtain a contradiction.

\item[(Case 2.2.2)]
For all $C>0,$ there exists $\tau \in (- \infty,0),$ and $\de \in (0,1),$ such that for all 
 $(x,t)\in (\Omega\timess ( -\infty, \tau)),$  $(x,t)$ is not a $C$-essential $\de$-necklike point.
Choose $C \leq {1 \on 16},$ 
and let $\tau,\de$ be the $\tau,\de$ from the statement at the beginning of this case.
Set $$G := |t|^{\ep   \on 2} {| \trfrriem|^2 \on \Sc^{2 -\ep }},$$
with $\ep \leq \eta(\de):= {\de \on 100 (3-\de)}$ (notice that 
this function is well defined,as $\Sc >0$ everywhere).
Then, as Chow and Knopf show in \cite{ChKn} (see the proof of theorem 9.19 there)
\begin{eqnarray}
 \partt G \leq \lap G + 2 {(1 - \ep) \on \Sc } \langle \grad G, \grad \Sc \rangle - {\ep \on 2 |t| } G,\label{Geqn} \end{eqnarray} 
for all $t \leq \tau.$
Let us examine $G$  a little more carefully. For fixed $t <0$ and a fixed $x_0$ we have the estimate
\begin{eqnarray}
\lim_{d(x,x_0,t) \to \infty} G(x,t) &=& 
\lim_{d(x,x_0,t) \to \infty}|t|^{\ep   \on 2} {| \trfrriem(x,t) |^2 \on \Sc^2(x,t) } \Sc^{\ep}(x,t) \cr
&\leq&  |t|^{\ep   \on 2} c(n)\lim_{d(x,x_0,t) \to \infty} \Sc^{\ep}(x,t) \cr
&=&0 \label{GD}
\end{eqnarray}
in view of the fact that the asymptotic scalar curvature ratio is less than infinity.
Also,  as Chow and Knopf point out, we have
\begin{eqnarray}
 G = |t|^{ \ep} \Sc^{\ep} {| \trfrriem |^2  \on  \Sc^2} {1 \on |t|^{\ep \on 2} }
 \leq {B^{\ep} c(n)  \on  |t|^{\ep \on 2} }, \label{Geqn2}
\end{eqnarray}
in view of the fact that 
$ B:= \sup_{\Omega\timess ( -\infty, \omega]}  |t| |\Sc(t)| < \infty,$
and hence 
\begin{eqnarray}
\lim_{t \to -\infty} \sup_{x \in M} G(x,t) = 0. \label{GT}
\end{eqnarray}

Let $\tau' < \tau -2$ be a constant with $\sup_\Omega G(\cdot,t) < \ep_0$ for all $t \leq \tau'.$
We know that
 \begin{eqnarray}
\sup_{M \timess (-\infty,0]} |\Riem|  \leq c(n)  \label{timeR} 
\end{eqnarray}
and without loss of generality
\begin{eqnarray}
 \sup_{M \timess [\tau', \tau] } |\grad \Riem|^2 + |\grad^2 \Riem|^2 \leq c(n)  \label{spaceR} 
\end{eqnarray}
in view of the interior gradient estimates of Shi (see  \cite{Ha95}, chapter 13).
We also know that for given $\ep_1 >0 $ and $s \in  [\tau', \tau] $ there exists an $r(s,\ep_1) >0$ such that
\begin{eqnarray}
 \sup_{ \{x \in M: d^2(x,x_0,s) \geq r \} } |\Riem|(x,s)  \leq \ep_1,\label{spaceG}
\end{eqnarray}
in view of the fact that the asymptotic scalar curvature ratio is finite.
Hence, for all $\ep_2 >0$  there exists a $ \de > 0,$ such that
$$ \sup_{x \in M, t \in (s,s + \de) : d^2(x,x_0,s) \geq r} |\Riem|(x,t)  \leq \ep_1 + \ep_2,$$
in view of \eqref{timeR} and  \eqref{spaceR} and the evolution equation for  $|\Riem|^2$.
In particular if   $\sup_{M} G(\cdot,s) < \ep_0,$ then
 $\sup_{M \timess (s,s + \de) } G(\cdot,t) < \ep_0,$ for small enough $\de$
(outside of a fixed large compact set $K$, $G< \ep_0$ for all $t \in (s,s+ \de)$ and inside $K$ we use the fact that $G$ is smooth). 
That is, the set $$Z:= \{r: \sup_{\Omega} G(\cdot,t) < \ep_0 \forall t \in [\tau', r) \}$$ is open.
Hence either 
  $$\sup_\Omega G(\cdot,t) < \ep_0$$ for all $t \in [\tau' ,\tau),$
or there is a first time $t_0 \in (\tau',\tau)$ such that 
$\sup_\Omega G(\cdot,t_0) = \ep_0.$
In the second  case, we see (using  equation \eqref{spaceG} with $s = t_0$ ) 
that there must also be a point $x_0 \in M$ such that  $ G(x_0,t_0) = \ep_0$.
But this contradicts the maximum principle in view of \eqref{Geqn}.

This means that $$\sup_\Omega G(\cdot,t) < \ep_0,$$ for all $t \in (-\infty,\tau),$
and hence, since $\ep_0$ was arbitrary,
$$ G \equiv 0.$$
Hence $\Omega = S^3 / \Gamma,$ which is a contradiction to the fact that
$\Omega$ is non-compact.

\end{itemize}

\end{proof}

\section{An application of the proof of Lemma \ref{globaltimelemma}. }

In certain cases, the proof of Lemma \ref{globaltimelemma} is applicable even if $M$ is non-compact.
For example, the theorem below is proved similarly to Lemma  \ref{globaltimelemma}.
This theorem was initially proved (using other methods) by Perelman Proposition 11.4 \cite{Pe1}.
\begin{theo}\label{globaltimelemmaapp}
Let $(\Omega^3,g(t))_{t \in (-\infty,0]}$ be an ancient non-compact complete solution to Ricci flow, with
(for some fixed origin  $o \in M$)
 \begin{eqnarray}
&& \sec \geq 0 \cr \label{a} \cr
&& \sup_\Omega |\Riem(g(t))|  < \infty \ \ \forall  t \in (-\infty,0) \label{b} \cr
&& \curlV(\tau):= \lim_{r \to \infty} {\vol(B_r(o,\tau)) \on r^n}   \geq \curlV_0 >0. \label{c} 
\end{eqnarray}
for some time $\tau,$  $\tau \in (-\infty,0).$ 
Then $(\Omega^3,g(t))$ is flat for all $t \in (-\infty,0)$.
\end{theo}

\begin{rema}
The limit in the statement of the theorem exists in view of the fact that
${\vol(B_r(o,\tau)) \on r^n}$ is non-increasing as $r$ increases (in view of the Bishop-Gromov comparison principle).
\end{rema}

\begin{proof}
Assume that the asymptotic scalar curvature ratio 
 $A = \limsup_{s \to \infty} \Sc s^2  = \infty$
(this is a constant independent of time).

Notice that for this solution, and any scaling of this solution which has bounded curvature
in a ball of radius one around some origin
$o'$, we have a uniform bound on the injectivity radius from below
at $o'$ (and time $s$),
in view of \ref{c} and \cite{CGT}: we have the estimate
$$ {\vol(B_r(o',s)) \on r^n}   \geq \curlV_0 >0$$
for all $r>0$ in view of \ref{c} and the Bishop Gromov volume comparison principle.
Furthermore $ {\vol(B_r(o',s)) \on r^n}   \leq \omega_n  $
trivially using the Bishop Gromov volume comparison principle. We may then apply the result 
of \cite{CGT} to obtain our estimate for the bound on the injectivity radius, exactly as we did in the argument of
Lemma \ref{globaltimelemma}.
Also, the estimates 
\begin{eqnarray}
\omega_n \geq  {\vol(B_r(o,s)) \on r^n}  \geq \curlV_0 >0 \forall r \geq 0 \label{skaly}
\end{eqnarray}
remain valid under scaling
(as the inequality is scale invariant).
Hence, we obtain a uniform bound from below on the injectivity radius estimate at $o'$, for 
any scaling of this solution which has bounded curvature by some fixed constant $c$ on a ball of radius one around $o'$

Translate in time so that $\tau = 0.$
Then we use the dimension-reduction argument of Hamilton (see Lemma 22.2 of \cite{Ha95} and the argument directly after 
the proof of Lemma 22.2) to obtain a new solution
 $ (N\timess \R , \ga\timess dx^2) $ or a quotient thereof by
group of isometries (without loss of generality, we may assume that we don't have a quotient, as explained in case 1.1 
of the proof of Lemma 
\ref{globaltimelemma}).
Notice that the dimension-reduction argument of Hamilton is applicable, in view
of the bounds from below on the injectivity radius at the centres of the balls occurring in the argument 
(due to the argument at the beginning of this theorem).
Also \ref{skaly} remains true for the resulting solution, as \ref{skaly} is scale invariant.
Without loss of generality the solution is defined on 
$ (N\timess \R , \ga\timess d\si^2)$ for $t \in (-\infty,\omega]$ for some $\omega >0,$ in view of the short time existence 
result of Shi, \cite{Sh}.

Assume that $\Sc(0,o) \neq 0$ on $N \timess (-\infty,\omega)$.
Then, see Lemma 26.2 in \cite{Ha95},
we have $$A = \limsup_{s \to \infty} \Sc s^2  < \infty $$ is a constant independent of $t \in (-\infty,\omega)$
on $N.$
This means that $$ \curlV(t) = \lim_{r \to \infty} {\vol(B_r(o,t)) \on r^n} $$ is a constant on $N$ independent of time, 
and in particular $$  \omega_n \geq {\vol(B_r(o,t)) \on r^n}  \geq \curlV_0 >0 \forall r \geq 0 \ \ \forall t \in (-\infty,\omega)$$ 
 (see theorem 18.3 in \cite{Ha95} ).

We then consider the following two cases:  
\begin{itemize}
\item[(case 1)]  $\sup_{N\timess ( -\infty, \omega]} |t| |\Sc(t)| = \infty, $
\item[(case 2)]  $ \sup_{N \timess ( -\infty, \omega]}  |t| |\Sc(t)| < \infty.$
\end{itemize}
exactly as in then proof of Lemma \ref{globaltimelemma}.
Both cases lead to a contradiction.

In the case that 
$A = \limsup_{s \to \infty} \Sc s^2 < \infty $ then we also know that 
$$ \curlV(t) = := \lim_{r \to \infty} {\vol(B_r(o,t)) \on r^n} $$ is a constant on $\Omega$ 
independent of time, 
and in particular$$   \omega_n \geq {\vol(B_r(o,t)) \on r^n}  \geq \curlV_0 >0 \forall r \geq 0 \ \ \forall t \in (-\infty,0).$$ 
Translate in time so that the solution is defined on $(-\infty,\omega),$ $\omega >0.$  
We then consider the following two cases.  
\begin{itemize}
\item[(case 1)]  $\sup_{\Omega\timess ( -\infty, 0])} |t| |\Sc(t)| = \infty, $
\item[(case 2)]  $ \sup_{\Omega \timess ( -\infty, 0]}  |t| |\Sc(t)| < \infty.$
\end{itemize}
exactly as in then proof of Lemma \ref{globaltimelemma}.
Both cases lead to a contradiction.

\end{proof}

\section{Bounds on the Ricci curvature from below under Ricci flow}
\setcounter{equation}{0}                                        

We prove quantitative estimates that tell us how quickly the Ricci curvature can decrease, if we assume
at time zero that the Ricci curvature is not too negative.
Both lemmas may be read independently of the rest of the results in this paper.

The first lemma is suited to the case that we have a sequence of solutions to Ricci flow
 $(M_i,{^i g}(t) )_{t \in [0,T)}$ whose initial data satisfies
\begin{eqnarray}
&&\Ricci( {^i g}(0))  \geq -\ep_i \Sc( {^i g}(0)) {^i g}(0),  
\end{eqnarray}
where $\ep_i \to 0$ as $i \to \infty.$
One application of this lemma is: if a subsequence of subsets $(\Omega_i,{^i g}(t) ), t \in {[0,T)}$ ($\Omega_i$ open) converges
(in the sense of Hamilton, see \cite{Ha95a})  to a smooth solution
$ (\Omega,{g}(t) ),t \in {(0,T)},$ then the lemma tells us that the Ricci curvature
of $ (\Omega,g(t))$ is non-negative for all $ t \in (0,T).$
This is very general, but does require that a limit solution exist. 

The second lemma is suited to the case that we have a sequence of solutions to Ricci flow
 $(M_i,{^i g}(t) )_{[0,T)}$ whose initial data satisfies
\begin{eqnarray}
&& \Ricci( {^i g}(0))  \geq -\ep_i {^i g}(0),
\end{eqnarray}
where $\ep_i \to 0$ as $i \to \infty.$
Once again, one application of this lemma is: if a subsequence of subsets $(\Omega_i,{^i g}(t) ), t \in {(0,S)}$ converges
(in the sense of Hamilton, see \cite{Ha95a}) to a smooth solution
$ (\Omega,{g}(t) ), t \in {(0,S)},$ then the lemma tells us that the Ricci curvature
of $ (\Omega,g(t))$ is non-negative for all $ t \in (0,S).$
Another useful application of the second lemma is:
if a solution $(M, g(t)), t \in { [0,T)}$ satisfies
\begin{eqnarray}
&&|\Riem( {g} )|  \leq {c_0 \on t}  \cr
&& \Ricci( { g}(0))  \geq -\ep { g}(0) 
\end{eqnarray}
then for a well controlled time interval the solution satisfies
$$\Ricci( {g}) \geq  - c_0 \ep  {g}.$$
As we saw in Lemma \ref{globaltimelemma}, such a bound is relevant to the question of existence
of solutions to the Ricci flow. 
We apply this lemma in the main Theorem \ref{realmaintheo} and the Application
\ref{appli}.

\begin{lemm}\label{seclemma}
Let $g_0$ be a smooth metric on a 3-dimensional manifold $M^3$ which satisfies
 \begin{eqnarray}
&& \Ricci(g_0) \geq - { \ep_0 \on 4}  \Sc g_0  \ \  ( \sec(g_0)  \geq - \Sc {\ep_0 \on 4}g_0 ) 
\end{eqnarray}

for some $0 < \ep_0 < {1 \on 100}$,  and 
let $(M,g(\cdot,t))_{t \in [0,T)}$ be a solution to Ricci flow with $g(0) = g_0(\cdot)$.
Then  \begin{eqnarray*}
\Ricci(g(t)) & \geq -\ep_0(1 + 4t) g(t) - \ep_0(1 + 4t) \Sc(g(t)) g(t), \fall t \in [0,T) \cap [0,{1 \on 8 }) \cr
 (\sec(g(t)) & \geq  -\ep_0({1 \on 2} + t) g(t) - \ep_0 ({1 \on 2} + t)\Sc(g(t)) g(t) , \fall t \in [0,T)\cap [0,{1 \on 8 }) )\cr
\end{eqnarray*}

\end{lemm}

\begin{proof}
Define $\ep = \ep(t) = \ep_0(1 + 4t),$
and the tensor $L(t)$ by 
$$ L_{ij} := \Ricci_{ij} + \ep \Sc g_{ij} + \ep g_{ij}.$$
We shall often write $\ep$ for $\ep(t)$ (not to be confused with $\ep_0$).
Notice that 
$\ep_0 < \ep(t) \leq 2\ep_0,$ for all $t \in  [0,{1 \on 8}):$ we will use this freely.
Then ${L_i}^j = ({R_i}^j + \ep R {\de_i}^j +  \ep{\de_i}^j),$ and  
 \begin{eqnarray*}
(\partt L)_{ij}  & = &
    (\partt L_i^l) g_{jl} - 2L_i^lR_{jl} \\
 & = & g_{jl} \Big( \partt ( R_{ik}g^{kl})  +  \ep \partt R {\de_i}^l + 4 \ep_0  R  {\de_i}^l +  4 \ep_0  {\de_i}^l\Big) - 2L_i^lR_{jl} \\
 &=& g_{jl}\partt ( R_{ik}g^{kl})   +   \ep g_{ij}\partt R + 4 \ep_0  R g_{ij} + 4 \ep_0 g_{ij}  - 2L_i^lR_{jl} \\
& = &  g_{jl}\Big( {(\lap \Ricci)_i}^l - {Q_i}^l + 
2R_{ik}R_{sm}g^{km}g^{ls} \Big)  \\
&& + \ep g_{ij}\Big( \lap R + 2|\Ricci|^2\Big)  +   4 \ep_0  R g_{ij}   + 4 \ep_0  g_{ij} - 2L_i^lR_{jl} \\
             &=& (\lap L)_{ij} - Q_{ij} + 2R_{ik}R_{jm}g^{km} +  2\ep|\Ricci|^2 g_{ij}\\
  && +4 \ep_0 R g_{ij} 
                  + 4 \ep_0  g_{ij} - 2L_i^lR_{jl}, \\
\end{eqnarray*}
where $Q$ is the tensor 
\begin{eqnarray}
&&Q_{ij} := 6S_{ij} - 3RR_{ij} + (R^2 -2S)g_{ij},\cr
&&S_{ij} := g^{kl}R_{ik}R_{jl} \label{Qeqn}
\end{eqnarray}
(see \cite{Ha82}, theorem 8.4)
Clearly $L_{ij}(0) > 0$. Define $N_{ij}$ by
$$N_{ij} := - Q_{ij} + 2R_{im}R_{sj}g^{ms} + 2 \ep |\Ricci|^2 g_{ij}+ 
4\ep_0 R g_{ij} + 4 \ep_0 g_{ij} - 2L_i^lR_{jl}.$$
We argue as in the proof of Hamilton's maximum principle,
Theorem 9.1, \cite{Ha82}.

We claim that $L_{ij}(g(t)) \geq 0.$
Assume there exist a first time and point $(p_0,t_0)$ and a direction $w_{p_0}$
for which $L(w,w)(g(t))(p_0,t_0) = 0.$  
Choose coordinates about $p_0$ so that at $(p_0,t_0)$ they are orthonormal, and 
so that $\Ricci$ is diagonal at $(p_0,t_0)$.
Clearly $L$ is then also diagonal at $(p_0,t_0)$.
W.l.o.g. 
\begin{eqnarray}
&&R_{11} = \la \cr
&&R_{22} = \mu\cr
&&R_{33}= \nu 
\end{eqnarray}
and
$\la \leq \mu \leq \nu,$
and so $$L_{11} =  \la +   \ep(t_0) R + \ep(t_0)\leq L_{22} \leq L_{33},$$
and so $L_{11} = 0,$ (otherwise $L(p_0,t_0) >0$ : a contradiction).
In particular, 
\begin{eqnarray}
N_{11}(p_0,t_0) =  
(\mu - \nu)^2 + \la(\mu + \nu) +  2\ep \la^2 + 2\ep\mu^2 + 2\ep\nu^2 + 
4\ep_0 R +  4 \ep_0, \label{N}
\end{eqnarray}
in view of the definition of $Q$ 
(see \cite{Ha82} Corollary 8.2, Theorems 8.3,8.4)
and the fact that $L_{11}=0$. Also, 
$L_{11} = 0 \Rightarrow \la = -\ep R - \ep$ 
at $(p_0,t_0)$,  and so, substituting this into \eqref{N}, we get

\begin{eqnarray*}
 N_{11}(p_0,t_0) & = & (u-v)^2 + (-\ep R - \ep)(\mu + \nu)
+ 2 \ep (\la^2 + \mu^2 + \nu^2)
+4 \ep_0 R + 4\ep_0 \\ 
&\geq&  \ep (   -( \la + \mu + \nu)(\mu + \nu) + 2 \la^2 + 2\mu
^2 + 2\nu^2)
+ 4 \ep_0 R + 4 \ep_0  - \ep(\mu + \nu)  \\
&=& \ep (   -( \la + \mu + \nu)(\mu + \nu) + 2 \la^2 + 2\mu^2 + 2\nu^2)
+ 4 \ep_0 R + 4 \ep_0  - \ep R + \ep \la \\
& \geq &  \ep (\la -\la \mu - \la \nu + \mu^2+ \nu^2 + 2\la^2 - 2 \mu \nu) + 4 \ep_0  R  + 4\ep_0 - \ep R. \cr
\end{eqnarray*}
To show $N_{11} > 0,$ we consider a number of cases.

\begin{itemize}
\item{ \bf Case 1}.
 $\la \geq 0.$
This combined with $L_{11} = 0$ implies that
$R < 0.$ A contradiction to the fact that $\la \geq 0$ and $\la$ is the smallest eigenvalue of $\Ricci$.
 \item{ \bf Case 2}.  $\la \leq 0, R \geq 0.$
 This implies $\nu \geq 0$ and hence
$$N_{11} \geq \ep (\la  -\la \mu + \mu^2+ \nu^2 + 2\la^2 - 2 \mu \nu)+  4 \ep_0,$$ 
in view of the fact that $\ep R \leq 2\ep_0 R.$
In the case $\mu \geq 0$ we obtain 
$$N_{11} \geq  \ep(\la + \mu^2+ \nu^2 + 2\la^2 - 2 \mu \nu) + 4 \ep_0 \geq -\ep + 4 \ep_0>0,  $$ 
after an application of Young's inequality, and similarly 
in the case  $\mu \leq 0$ we get 
$$ N_{11} \geq  \ep (\la  -\la \mu + \mu^2+ \nu^2 + 2\la^2) + 4 \ep_0 > 0.$$
\item{ \bf Case 3}.
 $\la \leq 0, R \leq 0.$
We know that $R(g_0) \geq - 3\ep_0$ will be preserved by Ricci flow,
and hence
$ 0 \geq R(g(t)) \geq - 3\ep_0.$
We break case 3 up into three sub-cases (3.1,3.2,3.3).
\begin{itemize}
\item{\bf Case 3.1} $\mu,\nu \leq 0.$
This with $R \geq - 3\ep_0$ implies that $| \la|,|\mu|, |\nu| \leq 3\ep_0$ and hence
 $$N_{11} \geq - 3\ep \ep_0 - 36 \ep \ep^2_0 -  12 \ep^2_0  + 4 \ep_0 \geq  -100 \ep^2_0 + 4 \ep_0 >0,$$ 
since $0< \ep_0 < {1 \on 100}, \ep < 2\ep_0 <1.$
\item{\bf Case 3.2} $\mu \leq 0 ,\nu \geq 0.$
Implies 
$$ N_{11} \geq  \ep (\la  -\la \mu+\mu^2+ \nu^2 + 2\la^2)  -  12 \ep^2_0 + 4 \ep_0 >0,$$
in view of Young's inequality, $\ep_0 \leq  {1 \on 100},$ and $0 < \ep < 2\ep_0.$
\item{\bf Case 3.3} $\mu \geq 0  (\Rightarrow \nu \geq 0).$
Then, similarly,  
$$N_{11} \geq  \ep ( \la +  \mu^2 + \nu^2 + 2\la^2 - 2 \mu \nu) -  12 \ep^2_0  + 4 \ep_0  
> 0.
$$
\end{itemize}

\end{itemize}
So in all cases $N_{11} >0.$
The rest of the proof is standard (see  \cite{Ha82} Theorem9.1):
extend $w(p_0,t_0)= \parti{}{x^1}(p_0,t_0)$ in space to a vector field $w(\cdot)$ in a small
neighbourhood of $p_0$ so that $^{g(t_0)}\grad w(\cdot)(p_0,t_0) = 0,$ and let
$w(\cdot,t) = w(\cdot).$
Then
$$0 \geq  (\partt L(w,w))(p_0,t_0) \geq  (\lap L(w,w))(p_0,t_0) + N(w,w) > 0 ,$$
which is a contradiction.
 
The case for the sectional curvatures is similar:
from \cite{Ha86}, Sec. 5, we know that the reaction equations
for the curvature operator are
\begin{eqnarray*}
\partt \al &=& \al^2 + \be\ga \\
\partt \be &=& \be^2 + \al \ga \\
\partt \ga &=& \ga^2 + \al \be.
\end{eqnarray*}
Note that: \begin{eqnarray}
R &=& \al + \be + \ga, \cr
|\Ricci |^2 &=& ( ({\al + \be \on 2})^2 + ({\al + \ga \on 2} )^2 + ({\be + \ga \on 2} )^2) \cr
&=&   {1 \on  2}( \al^2 + \be^2 + \ga^2 + \al \be + \al \ga + \be \ga). \label {rseceqn}
\end{eqnarray}

Similar to the Ricci case, we examine the function
$ \al + \ep(t) R  + \ep(t) $ where $\al \leq \be \leq \ga$ are eigenvalues of the curvature operator, and
$\ep(t) = \ep_0 ({1\on 2} + t).$ In order to make the following inequalities more
readable, we write
$\ep$ in place of $\ep(t)$: that is, $\ep =  \ep_0 ({1\on 2} + t)$.
\begin{eqnarray*}
\partt (\al + \ep R +  \ep) &=& \ep_0 + \ep_0 R + \al^2 + \be \ga + 2 \ep|\Ricci|^2 \\
                            &=& \ep_0 + \ep_0 R+ \al^2 + \be \ga 
                          + \ep( \al^2 + \be^2 + \ga^2 + \al \be + \al \ga + \be \ga),
\end{eqnarray*}
and so in the case that 
$\be, \ga \geq 0 ,$ or $\be, \ga \leq 0$, 
$\partt (\al + \ep R +  \ep) \geq \ep_0 (1 + R) >0.$ 
So assume that $\al \leq \be \leq 0,$ and $\ga \geq 0$.
Combining these inequalities with $\ep(t) \leq \ep_0 $, we see that
\begin{eqnarray*}
\partt (\al + \ep R + \ep )  &  \geq & \ep_0 +  \ep_0 R+ \al^2 +  \al \ga + \ep ( \al^2 + \be^2 + \ga^2  + \al \be+ \al \ga + \be \ga) \cr
                      & = & \ep_0 +  \ep_0 R + \al^2 + (\al + \ep R + \ep) \ga  \cr
                       &&  - \ep R \ga - \ep \ga  
                              + \ep ( \al^2 + \be^2  + \ga^2  + \al \be + \al \ga + \be \ga )  )\cr
                       &  = & \ep_0 +  \ep_0 R + \al^2 +   (\al + \ep R + \ep) \ga  - \ep \ga  
                          +\ep (  \al^2 + \be^2  + \al \be),\cr
                       & \geq & \al^2 +  (\al + \ep R + \ep) \ga + \ep_0(1 + R - \ga) +  \ep (\al^2 + \be^2 ),
\end{eqnarray*}
which, using $\ep(t) \geq {\ep_0 \on 2}$, is 
\begin{eqnarray*}
                       & \geq &  \al^2 + (\al + \ep R + \ep) \ga   + \ep_0 ( 1 +\al + \be + {\al^2 \on 2} + {\be^2 \on 2} ),\cr
                       & \geq &  \al^2 + (\al + \ep R + \ep) \ga,
\end{eqnarray*}
in view of Young's inequality.
At a point where $\al  + \ep R  +  \ep = 0,$
the last sum is strictly bigger than zero ( if $\al = 0$, then, $R \geq 0,$ and hence 
$\al + \ep R + \ep  \geq \ep > 0:$ a contradiction).
Then we argue as above.
 
\end{proof}          

The above lemma shows us that if the Ricci curvature at time zero is bigger than $-\ep$ ($\ep$ small)
then the  Ricci curvature divided by the scalar curvature is at most $-c\ep$
at points where the scalar curvature is bigger than one (for a short well defined time interval).
It can of course happen that the Ricci curvature becomes very large and negative in a short time,
if the scalar curvature is very large and positive in a short time.

Now we prove an improved version of the above theorem, which allows for some scaling in time.
In particular, for the class of solutions where $|\Riem|t \leq c_0$ it tells us that: if
the  Ricci curvature at time zero is bigger than $-\ep$ ($\ep$ small) 
then the  Ricci curvature is at most $-c\ep$ for some short well defined time interval.

\begin{lemm}\label{seclemma2}
Let $g_0$ be a smooth metric on a 3-dimensional manifold $M^3$ which satisfies
 \begin{eqnarray}
&& \Ricci(g_0) \geq - { \ep_0 \on 4} g_0, \cr
&& (\sec(g_0)  \geq - { \ep_0 \on 4} g_0 )
\end{eqnarray}
for some $0 < \ep_0 < {1 \on 100}$,  and 
let $(M,g(\cdot,t))_{t \in [0,T)}$ be a solution to Ricci flow with $g(0) = g_0(\cdot).$
Then  \begin{eqnarray*}
\Ricci(g(t)) & \geq -\ep_0(1 + kt) g(t) - \ep_0(1 + kt) t \Sc(g(t)) g(t), \fall t \in [0,T) \cap [0,T') \cr
 (\sec(g(t)) & \geq  -\ep_0({1 \on 2} + kt) g(t) - \ep_0 ({1 \on 2} + kt) t \Sc(g(t)) g(t) , \fall t \in [0,T)\cap [0,T'))
\end{eqnarray*}
where $k = 100$ and $T'=T'(100) > 0$ is a universal constant.
\end{lemm}

\begin{proof}
The proof is similar to that above.
Define $\ep = \ep(t) = \ep_0(1 + kt),$
and the tensor $L(t)$ by 
$$L_{ij} := \Ricci_{ij} + \ep t  \Sc g_{ij} + \ep g_{ij}.$$
We shall often write $\ep$ for $\ep(t)$ (not to be confused with $\ep_0$).
Notice that 
$\ep_0 < \ep(t) \leq 2\ep_0,$ for all $t \in  [0,{1 \on k}):$ we will use this freely.
Then $${L_i}^j = ({R_i}^j + \ep t R {\de_i}^j +  \ep{\de_i}^j),$$ and  
 \begin{eqnarray*}
(\partt L)_{ij}  & = &
    (\partt L_i^l) g_{jl} - 2L_i^lR_{jl} \\
 & = & g_{jl} \Big( \partt ( R_{ik}g^{kl})  +  \ep \Sc {\de_i}^l +  \ep t \partt \Sc {\de_i}^l + k \ep_0 t \Sc  {\de_i}^l +  k \ep_0  {\de_i}^l\Big) - 2L_i^lR_{jl} \\
 &=& g_{jl}\partt ( R_{ik}g^{kl}) +  \ep \Sc g_{ij}  +  \ep t g_{ij}\partt \Sc + k \ep_0 t \Sc g_{ij} + k \ep_0 g_{ij}  - 2L_i^lR_{jl} \\
& = &  g_{jl}\Big( {(\lap \Ricci)_i}^l - {Q_i}^l + 
2R_{ik}R_{sm}g^{km}g^{ls} \Big)  + \ep \Sc g_{ij} \\
&& + \ep t g_{ij}\Big( \lap R + 2|\Ricci|^2\Big)    + k \ep_0 t \Sc g_{ij}   + k \ep_0  g_{ij} - 2L_i^lR_{jl} \\
             &=& (\lap L)_{ij} - Q_{ij} + 2R_{ik}R_{jm}g^{km} + \ep  \Sc g_{ij} +  2\ep t |\Ricci|^2 g_{ij}\\
  && + k \ep_0 t R g_{ij}                 + k \ep_0  g_{ij} - 2L_i^lR_{jl}, \\
\end{eqnarray*}
where $Q$ is the tensor defined in Equation \eqref{Qeqn}.
Clearly $L_{ij}(0) > 0$. Define $N_{ij}$ by
$$N_{ij} := - Q_{ij} + 2R_{ik}R_{jm}g^{km} + \ep  \Sc g_{ij} +  2\ep t |\Ricci|^2 g_{ij}
+ k \ep_0 t R g_{ij}                 + k \ep_0  g_{ij} - 2L_i^lR_{jl}.$$
We argue as in the proof of Hamilton's maximum principle,
Theorem 9.1, \cite{Ha82}.

We claim that $L_{ij}(g(t)) > 0$ for all $t \in [0,T).$
Assume there exist a first time and point $(p_0,t_0)$ and a direction $w_{p_0}$
for which $L(w,w)(g(t))(p_0,t_0) = 0.$  
Choose coordinates about $p_0$ so that at $(p_0,t_0)$ they are orthonormal, and 
so that $\Ricci$ is diagonal at $(p_0,t_0)$.
Clearly $L$ is then also diagonal at $(p_0,t_0)$.
W.l.o.g. 
\begin{eqnarray}
R_{11} &=& \la, \cr
R_{22} &=& \mu, \cr
R_{33} &=& \nu, 
\end{eqnarray}
and 
$$ \la \leq \mu \leq \nu,$$
and so $$L_{11} =  \la +   \ep(t_0) t_0 R + \ep(t_0)\leq L_{22} \leq L_{33},$$
and so $L_{11} = 0,$ (otherwise $L(p_0,t_0) >0$ : a contradiction).
In particular, 
\begin{eqnarray}
N_{11}(p_0,t_0) =&&  
(\mu - \nu)^2 + \la(\mu + \nu) +  2\ep t \la^2 + 2\ep t\mu^2 + 2\ep t \nu^2 \cr
&&+
 \ep  \Sc g_{ij} + k \ep_0 t R g_{ij} + k \ep_0  g_{ij}
\end{eqnarray}
in view of the definition of $Q$ 
(see \cite{Ha82} Corollary 8.2, Theorems 8.3,8.4)
and the fact that $L_{11}=0$. We will show that $N_{11}(p_0,t_0) >0.$
$L_{11} = 0 \Rightarrow \la = -\ep t_0  R - \ep$ 
at $(p_0,t_0)$,  and so, substituting this into \eqref{N}, we get

\begin{eqnarray*}
 N_{11}(p_0,t_0) & = & (u-v)^2 + (-\ep t_0 R - \ep)(\mu + \nu)
+ 2 \ep t_0 (\la^2 + \mu^2 + \nu^2) \cr
&&+
\ep  \Sc  + k \ep_0 t \Sc g_{ij} + k \ep_0 
 \cr 
&\geq&  \ep t_0 (   -( \la + \mu + \nu)(\mu + \nu) + 2 \la^2 + 2\mu
^2 + 2\nu^2) - \ep(\mu + \nu)  \cr
 && + \ep  \Sc  + k \ep_0 t_0 \Sc  + k \ep_0 
\cr
&=& \ep t_0 (   -( \la + \mu + \nu)(\mu + \nu) + 2 \la^2 + 2\mu^2 + 2\nu^2) \cr
&& + ( (-\ep^2 t_0 + k \ep_0 t_0 ) \Sc - \ep^2
+ k \ep_0)  \\
 &\geq&  
 \ep t_0 ( - \la \mu - \la \nu + \mu^2 + \nu^2 + 2 \la^2 - 2 \mu \nu) \cr
&& +  ( (-\ep^2 t_0 + k \ep_0 t_0 ) \Sc - \ep^2
+ k \ep_0)
\end{eqnarray*}
where here we have used once again that  
$$\la(x_0,t_0) = -\ep(t_0) t_0  R(x_0,t_0) - \ep(t_0).$$
If $\Sc(x_0,t_0) \leq 0,$ then using the fact that $\Sc \geq - 3 \ep_0$
is preserved by the flow,
we see that $$\ (-\ep^2(t_0)t_0 + k \ep_0(t_0) t_0 ) \Sc(x_0,t_0) - \ep^2
+ k \ep_0  \geq {k \on 2} \ep_0.$$
Furthermore,
\begin{itemize}
\item{[i]} $\la = -\ep R -\ep \leq \ep$ (since $R \geq -3\ep_0$) and  $\la = -\ep R -\ep \geq   -\ep $, that is
$|\la| \leq \ep$.
\item{[ii]} 
Similarly $|\mu +\nu| = |R -\la| \leq 4 \ep.$  
\end{itemize}
Hence
$$\ep t_0 ( - \la (\mu + \nu) + \mu^2 + \nu^2 + 2 \la^2 - 2 \mu \nu) \geq - 50 \ep^2_0,$$ and so
$N_{11}(p_0,t_0)  > 0.$
Hence we must only consider the case $R(p_0,t_0) \geq 0.$

\begin{itemize}
\item{ \bf Case 1}.
 $\la \geq 0.$
This combined with $L_{11} = 0$ implies that
$R(p_0,t_0) < 0.$ A contradiction.
\item{ \bf Case 2}.
$\la \leq 0, \mu \geq 0, \nu \geq 0.$

In this case we trivially obtain 
$N_{11} > 0.$
\item{ \bf Case 3}.
 $\la \leq 0, \mu \leq 0,  \nu \geq 0.$
Implies $$N_{11} > \ep t_0 ( -\la \mu+\mu^2+ \nu^2 + 2\la^2)  \geq 0,$$
in view of Young's inequality.
\end{itemize}

So in all cases $N_{11} >0.$
The rest of the proof is standard (see  \cite{Ha82} Theorem9.1):
extend $w(p_0,t_0)= \parti{}{x^1}(p_0,t_0)$ in space to a vector field $w(\cdot)$ in a small
neighbourhood of $p_0$ so that $^{g(t_0)}\grad w(\cdot)(p_0,t_0) = 0,$ and let
$w(\cdot,t) = w(\cdot).$
Then
$$0 \geq  (\partt L(w,w))(p_0,t_0) \geq  (\lap L(w,w))(p_0,t_0) + N(w,w) > 0 ,$$
which is a contradiction.
 
The case for the sectional curvatures is similar:
from \cite{Ha86}, Sec. 5, we know that the reaction equations
for the curvature operator are
\begin{eqnarray*}
\partt \al &=& \al^2 + \be\ga \\
\partt \be &=& \be^2 + \al \ga \\
\partt \ga &=& \ga^2 + \al \be.
\end{eqnarray*}
In what follows, we use the formulae \eqref{rseceqn} freely.

Similar to the Ricci case, we examine the function
$ \al + \ep(t) t R  + \ep(t) $ where $\al \leq \be \leq \ga$ are eigenvalues of the curvature operator, and
$\ep(t) = \ep_0 ({1\on 2} + kt).$ In order to make the following inequalities more
readable, we write
$\ep$ in place of $\ep(t)$: 
that is, $\ep =  \ep_0 ({1\on 2} + kt)$.
We assume $t \leq {1 \on 2 k}$ so that $\ep_0 {1\on 2} \leq \ep(t) \leq \ep_0.$ 
\begin{eqnarray*}
\partt (\al + \ep t R +  \ep) &=& \ep R + k \ep_0 t R + k \ep_0 + \al^2 + \be \ga + 2 \ep t |\Ricci|^2 \cr
                           & =&\ep R + k \ep_0 t R + k \ep_0 + \al^2 + \be \ga \cr
                         && + \ep t ( \al^2 + \be^2 + \ga^2 + \al \be + \al \ga + \be \ga),
\end{eqnarray*}
and so in the case that 
$\be, \ga \geq 0 ,$ or $\be, \ga \leq 0$, 
\begin{eqnarray}
\partt (\al + \ep R +  \ep) &\geq&   \ep R + k \ep_0 t R + k \ep_0 \cr
 &\geq& - 3\ep^2_0 - 3 \ep^2_0 + k \ep_0 > 0
\end{eqnarray} 
So assume that $\al \leq \be \leq 0,$ and $\ga \geq 0$.
Combining these inequalities with $\ep(t) \leq \ep_0 $, we see that
\begin{eqnarray*}
\partt (\al + \ep t R + \ep )  &  \geq &  \ep R + k \ep_0 t R + k \ep_0 +  \al \ga \cr
            && + \ep t ( \al^2 + \be^2 + \ga^2  + \al \be+ \al \ga + \be \ga) \cr
                      & = &  \ep R + k \ep_0 t R + k \ep_0 + (\al + \ep t R + \ep) \ga  \cr
                       &&  - \ep t R  \ga - \ep \ga  
                              + \ep t ( \al^2 + \be^2  + \ga^2  + \al \be + \al \ga + \be \ga ) \cr
                        & = &   \ep R - \ep \ga + k \ep_0 t R + k \ep_0  +   (\al + \ep t R + \ep) \ga  \cr
                         &&      +\ep t (  \al^2 + \be^2  + \al \be)\cr
                        & = &   \ep ( \al + \be) + k \ep_0 t R + k \ep_0  +   (\al + \ep t R + \ep) \ga  \cr 
                          &&     + \ep t (  \al^2 + \be^2  + \al \be)\cr
                       & \geq &  ( 2\ep \al + k \ep_0 t R + k \ep_0  )  
                               +\ep t (  \al^2 + \be^2  + \al \be)\\            
\end{eqnarray*}
at a point where $\al  + \ep t R  + \ep = 0.$
Using  $\al  + \ep t R  + \ep = 0$ again, we get  
\begin{eqnarray*}
 2\ep \al + k \ep_0 t R + k \ep_0   &= & 2\ep( - \ep t R - \ep)   + k \ep_0 t R + k \ep_0 \cr
& = & Rt( - 2 \ep^2 + k \ep_0 )+ k \ep_0 -2 \ep^2 \cr
&>&  {k \on 2} \ep_0,
\end{eqnarray*}
 since $R \geq -3 \ep_0$ is preserved by the flow,
and $t \leq {1 \on k}.$
Hence
  \begin{eqnarray*}
\partt (\al + \ep t R + \ep )  &  \geq &  {k \on 2} \ep_0 + \ep t (  \al^2 + \be^2  + \al \be) >0,
\end{eqnarray*}
at a point where $\al  + \ep t R  + \ep = 0.$
Then we argue as above.
 
\end{proof}          

So although the Ricci curvature can become very large and negative under the Ricci
flow, it can only do so at a controlled rate.
In particular, as we mentioned before this lemma, 
if the curvature satisfies $|\Riem| t \leq c_0$ for all $t \in [0,T)$ (in addition to the initial conditions) then
$\Ricci \geq - c_1(c_0) \ep_0,$ is true on some well defined time interval
$[0,T')$ (in dimensions two and three).

\section{Bounding the diameter and volume in terms of the curvature}

\setcounter{equation}{0}

The results of this section hold for all dimensions.

\begin{lemm}\label{diamlemm}
Let $(M^n,g(t))_{t \in [0,T)}$ be a solution to Ricci flow
with \begin{eqnarray}
&&\Ricci(g(t)) \geq -c_0 \cr
&& | \Riem(g(t))| t \leq c_0 \cr
&& \diam(M,g_0) \leq d_0 \label{bobo2}
\end{eqnarray}
Then
\begin{eqnarray} 
d(p,q,0) - c_1(t,d_0,c_0,n) \geq  d(p,q,t) \geq  d(p,q,0) -c_2(n,c_0) \sqrt t \label{diogo}
\end{eqnarray}
for all $t \in [0,T)$, where
$$  c_1(t,d_0,c_0,n) \to 0$$ as $t \to 0.$

In particular if ${^i g}_0$ is a sequence of smooth metrics on manifolds $M_i$ with 
\begin{eqnarray}
&&\diam(M_i, {^i g}_0) \leq d_0 \cr
&&\GH( (M_i,d({^i g}_0)) , (X, d_X) ) \itoinfty 0 
\end{eqnarray}
and
$(M_i, {^i g}(t))_{t \in [0,T_i)} $ are solutions to Ricci flow with 
\begin{eqnarray}
&& {^i g}(0) = {^i g}_0 \cr
&&\sec( {^i g}(t)) \geq -c_0 \ \ \ ( \Ricci( {^i g}(t))  \geq -c_0 ) \cr
&& | \Riem ({^i g}(t)) | t \leq c_0 \ \forall t \in [0,T_i),  \label{bobo}
\end{eqnarray}
then
$$\GH((M_i,d({^i g}(t_i))), (X, d_X) ) \itoinfty 0 $$ for any sequence $t_i \in [0,T_i), i \in \N$ where $t_i \itoinfty 0$.
\end{lemm} 

\begin{proof}
The first inequality 
$$ d(p,q,t) \geq  d(p,q,0) -c_1(n,c_0) \sqrt t$$
is proved in \cite{Ha95}, theorem 17.2 (with a slight modification of the proof: see Appendix C).
The second inequality follows easily from
\cite{Ha95}, Lemma 17.3: see Appendix C.

The second statement is a consequence of the first result, and the triangle inequality which is valid for the Gromov
Hausdorff distance:
\begin{eqnarray}
&&\GH((M_i,d({^i g}(t_i))),   (X, d_X) )\cr
&& \leq  \GH( (M_i,d({^i g}(t_i))), (M_i,d({^i g}_0)) )  +  \GH ( (M_i,d({^i g}_0)),(X,d_X) ) \cr
&&\leq c(t_i)  +  \GH( (M_i,d({^i g}_0)),(X, d_X) ) \itoinfty 0. 
\end{eqnarray} 

Here we have used the characterisation of Gromov Hausdorff distance given in \ref{hausap},
and the fact that the identity map $I:(M_i,d({^i g}(t_i))) \to (M_i, d({^i g}_0) ),$
is an $c(t_i)$ -Hausdorff approximation, where $c(t) \to 0$ as $t \to 0$ : see Appendix A, Definition \ref{Haspro} and Lemma \ref{hausap}.
\end{proof}

\begin{coro}\label{volcorr}
Let  $(M^n,{g}(t))_{t \in [0,T)}$ be an arbitrary solution to Ricci flow ($g(0) = g_0$)
satisfying the conditions
\ref{bobo2}
and 
assume that there exists $v_0 >0$  such that 
\begin{eqnarray}
&&\vol(M,{g}_0) \geq v_0 >0.
\end{eqnarray}
Then there exists an $S = S(d_0,c_0,v_0,n) >0$ such that  $$\vol(M,{g}(t)) \geq {3v_0 \on 4} 
\ \ \forall \ \ t \in [0,T) \intersect [0,S) $$
\end{coro}

\begin{proof}
If this were not the case, then there exist solutions 
 $(M_i^n,{^i g}(t))_{t \in [0,T_i)}$ satisfying the stated conditions and there exist
$t_i \in [0,T_i),$ $t_i \itoinfty 0$
such that  
$\vol(M_i,{^i g}(t_i)) = {3v_0 \on 4}.$
But then $$\GH( (M_i,d({^i g}(t_i))) , (X, d_X)) \itoinfty 0 $$ from the lemma above.
According to \cite{BGP}, thm 10.8 for the case that  $\sec( {^i g}(t)) \geq -c_0$ 
(for the Ricci case we use theorem
5.4 of \cite{ChCo}  Cheeger-Colding )
we also have $$ v_0 \leq \vol(M_i,{^i g}_0) = \Haus^n(M_i,d({^i g}_0) ) \itoinfty \Haus^n(X,d_X)$$ which implies
$ \Haus^n(X,d_X) \geq v_0.$ 
Here $\Haus^n(X,d_X)$ is the $n$-dimensional Hausdorff mass of $X$ with respect to the metric $d_X$.
Similarly we have  
$$  {3v_0 \on 4} = \Haus^n(M_i,d({^i g}(t_i))) \itoinfty \Haus^n(X,d_X)$$
Implies $\Haus^n(X,d_X) =  {3v_0 \on 4}.$
A contradiction.
\end{proof}

\section{Non-collapsed compact three manifolds of almost non-negative curvature.}

\setcounter{equation}{0}

The results of this section are only valid for dimensions two and three-


\begin{theo}\label{realmaintheo}
Let $M$ be a closed three (or two) manifold satisfying
\begin{eqnarray*}
&&\diam(M,g_0) \leq d_0 \cr
&&\Ricci(g_0) \ \ \ (\sec(g_0) )  \geq  - \ep g_0 \cr
&&\vol(M,g_0) \geq v_0  >0,
\end{eqnarray*}
where $\ep \leq  {1 \on 10 c^2}$ and $c = c(v_0,d_0) \geq 1$ is the constant from Lemma 
\ref{globaltimelemma}
Then there exists an $S = S(d_0,v_0) >0$ and $K = K(d_0,v_0)$ such that the
maximal solution $(M,g(t))_{t \in [0,T)}$ to Ricci-flow satisfies
$T\geq S,$
and $$\sup_{M} |\Riem(g(t))| \leq {K \on t},$$ for all
$t \in (0,S)$.

\end{theo}

\begin{proof}

Let $[0,T')$ be the maximal time interval for which
\begin{eqnarray*}
\vol(M,g(t)) &>& {v_0 \on 2}, \\
\Ricci(g(t)) &\geq& -1 \\
\diam(g(t)) &\leq& 5d_0. \\
\end{eqnarray*}
If $T' \leq 1,$
then the diameter condition will not be violated as long as the other conditions are not violated
(as one easily sees by examining the evolution equation for distance under Ricci flow). 
So we assume w.l.o.g $T' \leq 1$ and the diameter condition
is not violated.
From Lemma \ref{globaltimelemma}, we know that there exists a $c = c(d_0,v_0)$ such that
$R(t) \leq {c \on t},$ for all $ t \in [0,T').$
Using Lemma \ref{seclemma2} we see that there exists a $T'' = T''(c)>0$ such that
$\Ricci \geq -{1 \on 2}$ for all $t \in [0,T'']\cap [0,T')$. 
So the Ricci curvature condition is not violated on $[0,T'']\cap [0,T')$.
Furthermore, in view of Corollary \ref{volcorr} there exists a $T'''= T'''(v_0,d_0,c)$,
such that $\vol(M,g(t)) > {3v_0 \on 4}$ for all $t \in [0,T'''] \cap [0,T''] \cap [0,T']$
Hence $T' \geq \min(T''(c), T'''(v_0,d_0)) > 0,$ as required.
The estimate for the curvature and the existence of $S$ then follow from Lemma \ref{globaltimelemma}.

\end{proof}
 
\begin{theo}\label{appli}
Let $(M_i, \upi g_0)$ be a sequence of closed three (or two) manifolds satisfying
\begin{eqnarray*}
&&\diam(M_i,\upi g_0) \leq d_0 \cr
&&\Ricci(\upi g_0) (\sec(\upi g_0) )  \geq  - \ep(i) {\upi g}_0 \cr
&&\vol(M,\upi g_0) \geq v_0 >0,
\end{eqnarray*}
where $\ep(i) \to 0,$ as
$i \to \infty.$

Then there exists an $S = S(v_0,d_0) >0$ and $K = K(v_0,d_0)$ such that the
maximal solutions $(M_i,\upig(t))_{t \in [0,T_i)}$ to Ricci-flow satisfy
$T_i \geq S,$
and $$\sup_{M} |\Riem(\upig(t))| \leq {K \on t},$$ for all
$t \in (0,S)$.
In particular the Hamilton limit solution 
$(M,g(t))_{t \in (0,S)} = \lim_{i \to \infty} (M_i,\upig(t))_{t \in (0,S)}$ (see \cite{Ha95a})
exists and satisfies
\begin{eqnarray}
&& \sup_{M} |\Riem(g(t))| \leq {K \on t} \cr
&& \Ricci(g(t)) \geq 0  \ \ (\sec(g(t)) \geq 0),
\end{eqnarray}
for all $t \in (0,S)$ and $(M,g(t))$ is closed. Furthermore 
\begin{eqnarray}
\GH( (M,d(g(t))) , (M, d_{\infty}) ) \to 0 
\end{eqnarray}
as $t \to 0$ where
$(M,d_\infty) = \lim_{i \to \infty} (M_i,d(\upi g_0) )$
(the Gromov-Hausdorff limit). Hence, if $M = M^3$, then $M^3$ is diffeomorphic to
to a quotient of one of $S^3$, $ S^2 \timess \R $ or $ \R^3$ by a finite group
of fixed point free isometries acting properly discontinuously. 
\end{theo}

\begin{proof}

We apply the previous theorem. Then notice that
lemma \ref{seclemma} (or lemma \ref{seclemma2}) implies that $\Ricci(g(t)) \geq 0$ ($\sec(g(t)) \geq 0$) for this
limit solution, for all $t \in (0,S).$
To prove that $\GH( (M,d(g(t))) , (M, d_\infty) ) \to 0 $
use the triangle inequality as in the proof of Lemma \ref{diamlemm}:
\begin{eqnarray}
&&\GH((M,d({g}(t))),   (M, d_\infty) ) \cr
 &&\leq \GH((M,d({g}(t))),(M_i,d({^i g}(t))) +   \GH((M_i,d({^i g}(t))),   (M, d_\infty) )\cr
 &&\leq \GH((M,d({g}(t))),(M_i,d({^i g}(t)))+  
\GH( (M_i,d({^i g}(t))), (M_i,d({^i g}_0)) ) \cr
&&  +  \GH ( (M_i,d({^i g}_0)),(M,d_\infty) ) \cr
&&\leq \GH((M,d({g}(t))),(M_i,d({^i g}(t))) + c(t)  +  \GH( (M_i,d({^i g}_0)),(M, d_\infty) ) \cr
&& \itoinfty c(t),\end{eqnarray} 
for all $t >0$, where $c(t) \to 0$ as $t \to 0:$
here we have used \ref{diogo}, and the characterisation of Gromov-Hausdorff distance given in
\ref{hausap} to obtain $c(t).$
\end{proof}


\appendix

\section{Gromov Hausdorff space and Alexandrov spaces}

\begin{defi}
Let $(Z,d)$ be a metric space, $p \in Z,$ $r >0$.
$$B_{r}(p):= \{ x \in Z: d(x,p) <r \}.$$
For two non-empty subsets $A,B \subset Z$
$$\dist(A,B) = \inf\{ d(a,b): a  \in A, b \in B \} .$$
$$B_{r}(A):= \{ x \in Z: \dist({x},A) < r \}.$$
\end{defi}

\begin{defi}
For subsets $X,Y \subset (Z,d)$
we define the Hausdorff distance between $X$ and $Y$ by
$$d_H(X,Y):= \inf\{ \ep >0 : X \subset B_{\ep}(Y) \ \mbox{and} \ Y \subset B_{\ep}(X) \}.$$
\end{defi}

Then, (see \cite{BuBuIv} Prop.7.3.3)
\begin{prop}
\begin{itemize}
\item $d_H$ is a semi-metric on $2^Z$ (the set of all subsets of $Z$),
\item $d_H(A,\bar A) = 0$ for all $A \subset Z$, where $\bar A$ is the closure of $A$ ( in
$(Z,d)$)
\item If $A$ and $B$ are closed subsets of $(Z,d)$ and $d_H(A,B) = 0$ then $A = B$.
\end{itemize}
\end{prop}

\begin{defi}
For a subset $X \subset Z,$ $(Z,d)$ a metric space, we define $d|_X$ to be the metric on $X$ defined by 
$$d|_X(a,b) = d(a,b).$$
\end{defi}

We then define the Gromov-Hausdorff distance between two abstract metric spaces $(X,d_X)$ and $(Y,d_Y)$
as follows:
\begin{defi}

$d_{GH}( (X,d_X), (Y,d_Y) )$ is the infimum over all $ r > 0$ such that there exists a metric space $(Z,d)$
and maps $f:X \to Z$, $X' := f(X)$, and $g: Y \to Z$, $Y' := g(Y)$ such that
$f:(X,d_X) \to (X',d|_{X'})$ and $g:(Y,d_Y) \to (Y',d_{Y'})$ are isometries and 
$d_H(X',Y') < r.$
\end{defi}

\begin{fact} 
$d_{GH}$ satisfies the triangle inequality,i.e.,
$$ d_{GH}((X_1,d_1),(X_3,d_3)) \leq d_{GH}((X_1,d_1),(X_2,d_2)) + d_{GH}((X_2,d_2),(X_3,d_3))$$
 for all metric spaces $(X_1,d_1),(X_2,d_2),(X_3,d_3).$
\end{fact}

\begin{proof}
See \cite{BuBuIv}  Prop.7.3.16.
\end{proof}

\begin{defi}
An $\nu$-Hausdorff approximation $f:X \to Y$ for metric spaces $(X,d_X)$ and $(Y,d_Y)$ is
a map which satisfies
\begin{eqnarray}
&& |d_Y(f(x),f(x')) - d_X(x,x') | \leq \nu \cr
&& B_{\nu}(f(X)) = Y
\end{eqnarray}
\end{defi}
 
\begin{defi}\label{Haspro}
$\Happrox( (X,d_X), (Y,d_Y) )$ is the infimum of $\nu$ such that there exists a 
$\nu$-Hausdorff approximation $f:X \to Y$.
\end{defi}

We prove the following simple well known lemma
\begin{lemm}\label{hausap}
$$\Happrox( (X,d_X), (Y,d_Y) ) \leq 2 d_{GH}( (X,d_X), (Y,d_Y) \leq 4\Happrox( (X,d_X), (Y,d_Y) ).$$
\end{lemm}

\begin{proof}
Let  $f:X \to Y$ be a $\nu$-Hausdorff approximation.
Define $Z = X \disjointunion Y,$ where $\disjointunion$ denotes disjoint union, and define a metric there
\begin{eqnarray}
&& d_Z(x,x') = d_X(x,x') \ \mbox{if} \  x,x' \in X, \cr
 && d_Z(y,y') = d_Y(y,y') \ \mbox{if} \  y,y' \in Y, \cr
 &&d_Z(x,y) =  d_Z(y,x) = d_Y(y,f(x)) + \nu\ \mbox{if} \  y \in Y, x \in X 
\end{eqnarray}
We check that this defines a metric:
\begin{itemize}
\item[{(i)}]
$d_Z(a,b)\geq 0$ follows from the definition.

\item[{(ii)}]
$d_Z(a,b) = d_Z(b,a)$ per definition.

\item[{(iii)}]
$d_Z(a,b) = 0$ implies $a,b \in X$ with  $d_X(a,b) = 0$ or
 $a,b \in Y$ with  $d_Y(a,b) = 0,$ and hence, in both cases
$a =b$. $a = b$ implies $d_Z(a,b) = 0$ trivially.
\item[{(iv)}]
Assume $x,x' \in X$ and $y \in Y.$
Then 
\begin{eqnarray}
d_Z(x,y) + d_Z(x',y) &=& d_Y(f(x),y) + d_Y(f(x'),y) + 2 \nu \cr
 &\geq&   d_Y(f(x),f(x'))  + 2 \nu \cr
& \geq & d_Y( x,x'),
\end{eqnarray}
in view of the fact that  $f:X \to Y$ is a $\nu$-Hausdorff approximation.
For $y,y' \in Y,$$x \in X$ we have
\begin{eqnarray}
d_Z(x,y) + d_Z(x,y') &=& d_Y(f(x),y) + d_Y(f(x),y') + 2 \nu \cr
 &\geq&   d_Y(y',y)  + 2 \nu \cr
& \geq & d_Y(y,y'). 
\end{eqnarray}
The other cases which can occur also trivially satsify the Triangle inequality.
\end{itemize}

The inclusion maps $i_X: X \to X \disjointunion Y,$  $i_Y: Y \to X \disjointunion Y,$ 
are then isometries onto their images (here  $X \disjointunion Y$ is the disjoint union of $X$ and $Y$).
Furthermore
for $y \in Y$ there exists an $x_y \in X$ such that
$d_Y(f(x_y),y) \leq \nu$ as $f:X \to Y$ is a $\nu$ approximation.
Hence $$d_Z(x_y,y) =  d_Y(y,f(x_y)) + \nu \leq 2 \nu.$$
That is, $$ Y\subset B_{2\nu} (X ).$$
For $x \in X$ let $y_x := f(x)$. Then
 $$d_Z(x,y_x) =  d_Y(f(x),y_x) + \nu =  \nu.$$
Hence
$$ X \subset B_{2\nu} (Y ).$$
Hence 
$$d_{GH}( (X,d_X), (Y,d_Y) ) \leq 2 \nu .$$

For the other direction, assume that $$d_{GH}((X,d_X),(Y,d_Y)) \leq \nu.$$

Let $f:(X,d_X) \to (Z,d_Z),$ and $g:(Y,d_Y) \to (Z,d_Z)$ be one to one maps which are isometries onto
 their images,
and so that $Y' \subset B_{2\nu}(X')$ and $ X' \subset B_{2\nu}(Y'),$ where
$Y' = g(Y)$,$X' = f(X).$
That is, for each $x' \in X'$ there exists a $p(x') \in Y'$ with
$d(x',p(x')) \leq 2\nu.$ Define a function $p:X' \to Y'$ so that 
$d(p(x'),x') \leq 2\nu.$
 Define $l: X \to Y$ by 
$l(x) =:g^{-1}( p(f(x)) ).$
Then 
\begin{eqnarray}
d(x_1,x_2) - d(l(x_1),l(x_2) ) &=&     d( f(x_1), f(x_2) ) - d(p(f(x_1),p(f(x_2)) \cr
&\leq&  d( f(x_1), p(f(x_1)) ) +  d(  p(f(x_1)), p(f(x_2)) ) \cr
&& + \
d(  p(f(x_2)), f(x_2) ) - d(p(f(x_1),p(f(x_2)) \cr
&=& d( f(x_1), p(f(x_1)) ) + d(  p(f(x_2)), f(x_2) ) \cr
&\leq& 4 \nu.
\end{eqnarray}

Furthermore, let $y \in Y.$ Then $g(y) \in Y'$ and so there exists
an $x' \in X'$ with $d(x',g(y)) \leq 2\nu.$
Let $x \in X$ be the unique element with $f(x) = x'$.
Then 
\begin{eqnarray}
d(l(x),y) &=& d( g^{-1}( p(f(x)) ), y) \cr
&=&  d( p(x') , g(y) ) \cr
&\leq& d( p(x'), x') + d(x',g(y)) \leq d( p(x'), x') +  2\nu \cr
&\leq& 4 \nu 
\end{eqnarray}
where the last line follows in view of the definition of the function $p$.
Hence $$ Y \subset B_{4\nu}(l(X)).$$

\end{proof}

Now we state the compactness result of Gromov.

\begin{prop}
$\curlM(n,k,d_0)$ is precompact in Gromov-Hausdorff space. 
\end{prop}

\begin{proof}
See \cite{BuBuIv} Remark 10.7.5.
\end{proof}

Clearly $\curlS(n,k,d_0) \subset \curlM(n,(n-1)k,d_0)$ and so it is
also precompact in Gromov-Hausdorff space.

In \cite{BGP}
(Theorem 10.8), the following fact about the convergence of Hausdorff measure was shown.

\begin{theo}
Let $(M_i,g_i) \in \curlS(n,k,d_0)$,$i \in \N$  be a sequence of smooth Riemannian manifolds with 
$\vol(M_i,g_i) \geq v_0  >0,$ for all $i \in \N$ and
$$(M_i,d(g_i)) \itoinfty (X,d_X)$$ in Gromov Hausdorff space. 
Then 
$$\vol(M_i,g_i) = \Haus_{i}(M_i) \itoinfty  \Haus(M),$$ where
$\Haus_i: M_i \to \R_0^+$ is $n$-dimensional Hausdorff measure with respect to $d(g_i)$ and 
$\Haus: X \to \R_0^+,$ is  $n$-dimensional Hausdorff measure with respect to $d_X$.
\end{theo}

\begin{proof}
See for example Theorem 10.10.10 in \cite{BuBuIv}.
\end{proof}

In \cite{ChCo} (Theorem 5.4) the same result was proved for 
$\curlM(n,k,d_0).$

\begin{theo}
Let $(M_i,g_i) \in \curlM(n,k,d_0)$,$i \in \N$  be a sequence of smooth Riemannian manifolds
with $\vol(M_i,g_i) \geq v_0  >0$ for all $i \in \N$ , and
$$(M_i,d(g_i)) \itoinfty (X,d_X)$$ in Gromov Hausdorff space. 
Then $$\vol(M_i,g_i) = \Haus_{i}(M_i) \itoinfty  \Haus(M),$$ where
$\Haus_i: M_i \to \R_0^+$ is $n$-dimensional Hausdorff measure with respect to $d(g_i)$ and 
$ \Haus: X \to \R_0^+,$ is  $n$-dimensional Hausdorff measure with respect to $d_X$.
\end{theo}

\begin{proof}
See Theorem 5.4 of \cite{ChCo}.
\end{proof}

The spaces that arise in this chapter, are obtained as limits of spaces whose curvature is bounded
from below.
A.D. Alexandrov studied such spaces extensively, and there is a large field of literature
which is souly concerned with such spaces. 
Here we give one possible definition of the class of spaces with curvature bounded from below
(see \cite{BGP} and \cite{BuBuIv} for further possible definitions and properties of such spaces).
\begin{defi}
The complete metric space $(M,d)$ is called an intrinsic metric space if for any
$x,y \in M$, $\de >0$ there is a finite sequence of points $z_0 = x, z_1, \ldots, z_k = y,$
such that
$d(z_i,z_{i+1}) \leq \de \forall  i \in \{ 0, \ldots, k-1 \}$ and
$$ \sum_{i = 0}^{k-1} d(z_i,z_{i+1}) \leq d(x,y) + \de.$$
\end{defi}
\begin{defi}
The length $L_d(\ga)$ of a continuous curve $\ga:[a,b] \to M$ is the 
supremum of the sums
$$ \sum_{i = 0}^{k-1} d(\ga(y_i), \ga(y_{i+1}) )$$
 over all partitions $a=y_0 < y_1 < \ldots y_k = b,$ of
$[a,b]$ (notice that this length could be infinite). 
\end{defi}
\begin{defi}
A geodesic is a continuous curve $\ga:[a,b] \to M$ whose
length is equal to $d(\ga(a),\ga(b))$ (that is, the distance between the endpoints of the curve).
\end{defi}

\begin{defi}
A collection of three points $p,q,r \in M$ and three
geodesics $pq,pr,qr$ is called a triangle in $M$ and denoted by
$\tri(p,q,r)$.
\end{defi}

\begin{fact}\label{trifact}
Let $(M_k,d_k)$ denote the complete simply connected two dimensional Riemannian manifold of sectional curvature $k.$
Let us fix the real number $k$ and let $\tri(p,q,r)$ be given.
For $k<0$ there exists a unique (up to an rigid motion) triangle
$\tri(\ti p, \ti q, \ti r)$ in the metric space $(M_k,d_k)$ with
$d(p,q) = d_k(\ti p,\ti q), d(p,r) = d_k(\ti p, \ti r), d(r,q) = d_k(r,q).$
For $k>0$ we require that the perimeter of  $\tri(p,q,r)$ be less than
${2 \pi \on \sqrt k}$ in order that the triangle exist.
\end{fact}

\begin{defi}
A complete, locally compact space $M$ with intrinsic metric $d$ is
called an Alexandrov space with curvature $\geq k$  if in some neighbourhood  $U_x$ of each point $x$,
for any triangle $\tri(p,q,r)$ with vertices in $U_x$ and any point $s$ on the geodesic
$qr$, the inequality $d(p,s) \geq d_k(\ti p, \ti s)$ is satisfied, where
$\tri(\ti p, \ti q,\ti r)$ is the triangle from \ref{trifact}, and $\ti s$ is the
point in $\ti q \ti r$ satisfying $ d(q,s) = d_k(\ti q, \ti s)$ and $d(r, s) = d_k(\ti r, \ti s)$.  
\end{defi}

Many reuslts which are valid for smooth Riemannian manifolds with curvature bounded from below by $k$
are also valid for Alexandrov spaces with curvature $\geq k.$ 
For example, Theorem 3.6 of \cite{BGP} says
that for $(X,d)$ an Alexandrov space with curvature $\geq k$ we have
$$\diam(X,d) \leq {\pi \on \sqrt k}.$$
For other properties of Alexandrov spaces with curvature $\geq k$ see \cite{BGP} or the 
book \cite{BuBuIv}.

\section{$C$-essential points and $\de$-like necks}

\begin{defi}
Let $(M,g(t))_{t \in (-\infty,T)}$, $ T \in \R \cup \{ \infty \},$  be a solution to Ricci flow.
We say that $(x,t) \in M \timess (-\infty,T) $ is a $C$-essential  point if
$$ |\Riem(x,t)| |t| \geq  C.$$
\end{defi}
\begin{defi}
We say that  $(x,t) \in M\timess   (-\infty,T) $ is a $\de$-necklike point
if there exists a unit 2-form $\theta$ at $(x,t)$ such that
$$ | \Riem - R ( \theta \otimes \theta ) | \leq \de |\Riem|.$$
\end{defi}

$\de$-necklike points often occur in the process of taking a limit around a sequence of times and points
which are becoming singular.
If $\de=0,$ then the inequality reads
$$ | \Riem(x,t) - R(x,t) ( \theta \otimes \theta ) | = 0.$$
In three dimensions this tells us that the manifold splits.
This can be seen with the help of some algebraic lemmas.

\begin{lemm}
Let $\omega \in \Omega^2(\R^3).$
Then it is possible to write 
$$\omega = X \land V,$$ for two orthogonal vectors $X$ and $V$.
\end{lemm}

\begin{rema}
Here we identify one forms with vectors using 
$$ adx^1 + bdx^2 + cdx^3 \equiv (a,b,c).$$
\end{rema}

\begin{proof}

Assume
\begin{eqnarray}
\omega &=& a dx^1 \land dx^2 + b dx^1 \land dx^3 + c  dx^2 \land dx^3  
\end{eqnarray}
Without loss of generality $b \neq 0$. Then, we may write:
\begin{eqnarray}
\omega &=&  (dx^1 + {c \on b}dx^2)\land (a dx^2 + b dx^3)
\end{eqnarray}
So $\omega = X \land Y.$
Now let $X,Z,W$ be an orthogonal basis all of length $|X|.$
Then $$Y = a_1 X + a_2 Z + a_3 W.$$
This implies
\begin{eqnarray}
\omega &&=  X \land (a_1 X + a_2 Z + a_3 W) \cr
       &&=  X \land  (a_2 Z + a_3 W)
\end{eqnarray}
as required ($V =  a_2 Z + a_3 W$).
\end{proof}

Hence we may write the $\theta$ occurring above as
$$\theta =  X \land V.$$
Hence
$$\Riem(x,t) = c X \land V \otimes  X \land V,$$ with
$$\{ X,V,Z \}$$ an orthonormal basis for $\R^3.$

The set $\{ X \land V,  X \land Z, V \land Z  \}$ then forms an orthonormal basis
and the curvature operator $\curlR$ can be written with respect to this basis as
$$ \left( \begin{array}{rrr}
           c & 0 & 0  \\
           0  &  0 & 0 \\
           0 &  0 & 0
  \end{array} \right)
$$
Hence the manifold splits (if the solution is complete with bounded curvature and non-negative curvature operator) in view
of the arguments in chapter 9 of \cite{Ha86}.

\section{Estimates on the distance function for Riemannian manifolds evolving by Ricci flow}

For completeness, we prove some results which are implied or proved in \cite{Ha95} 
and stated in \cite{CCCY} as editors note 24 from the same paper in that book.
The lemma we wish to prove is 
\begin{lemm}\label{diamlemm2}
Let $(M^n,g(t))_{t \in [0,T)}$ be a solution to Ricci flow
with \begin{eqnarray}
&&\Ricci(g(t)) \geq -c_0 \cr
&& | \Riem(g(t))| t \leq c_0 \cr
&& \diam(M,g_0) \leq d_0 
\end{eqnarray}
Then
\begin{eqnarray} 
d(p,q,0) - c_1(t,d_0,c_0,n) \geq  d(p,q,t) \geq  d(p,q,0) -c_1(n,c_0) \sqrt t
\end{eqnarray}
for all $t \in [0,T)$, where
$$  c_1(t,d_0,c_0,n) \to 0$$ as $t \to 0.$

\end{lemm} 

\begin{proof}
The first inequality 
$$ d(p,q,t) \geq  d(p,q,0) -c_1(n,c_0) \sqrt t$$
is proved in \cite{Ha95}, theorem 17.2 after making a slight modification of the proof.
If we examine the proof there (as  pointed out in \cite{CCCY} as editors note 24 of the same book),
we see that in fact that what is proved is:
$$d(P,Q,t) \geq d(P,Q,0) -C \int_0^t \sqrt{M(t)}$$
where $\sqrt{M(t)}$ is any integrable function which satisfies
$$\sup_M |\Riem(\cdot,t)| \leq M(t).$$
In particular, in our case we may set $$M(t) = {c \on t}$$ which then implies the first inequality.
The second inequality is also a simple consequence of results obtained in \cite{Ha95}.
Lemma 17.3 tells us that
$$ \partt d(P,Q,t) \leq -\inf_{\ga \in \Gamma} \int_{\ga} \Ricci(T,T)ds$$
where the $\inf$ is taken over the compact set $\Gamma$ of all geodesics from $P$ to $Q$ realising the distance
as a minimal length, $T$ is the unit vector field tangent to $\ga.$
Then in our case $\Ricci \geq -c_0$ implies
$$ \partt d(P,Q,t) \leq c_0 d(P,Q,t).$$
This implies that 
$$ d(P,Q,t) \leq  \exp^{c_0 t}  d(P,Q,0),$$ and as a consequence
$$\diam(M,g(t)) \leq d_0 \exp^{ct}.$$
Hence 
\begin{eqnarray}
 d(P,Q,t) \leq  \exp^{c_0 t}  d(P,Q,0)&& =  d(P,Q,0) + (\exp^{c_0 t}-1) d(P,Q,0) \cr
                                      && \leq   d(P,Q,0) +(\exp^{c_0 t}-1)d_0 \exp^{ct},
\end{eqnarray}
which implies the result.

\end{proof}

\section{Notation}
\setcounter{equation}{0}
 $\R^+$ is the set of positive real numbers.\hfill\break
$\R_0^+$ is the set of non-negative real numbers.\hfill\break
For a Riemannian manifold $(M,g)$ $(M, d(g))$ is the metric space induced
by $g$.
For a tensor $T$ on $M$, we write
$\upg|T|^2$ to represent the norm of $T$ with respect to the metric $g$ on $M$.
For example if $T$ is a  ${0 \choose 2}$ tensor, then
$${}^g|T|^2 =  g^{ij}g^{kl}T_{ik}T_{jl}.$$\hfill\break
$\gradh T$ refers to the covariant derivative with respect to $h$ of $T$.\hfill\break
$\uph\Riem$ or $\Riem(h)$ refers to the Riemannian curvature tensor with respect to $h$ on $M$.\hfill\break
$\uph \Ricci$ or $\Ricci(h)$ or $\uph R_{ij}$ refers to the 
Ricci curvature of $h$ on $M$.\hfill\break
$\uph \Sc$ or $\Sc(h)$ refers to the scalar curvature of $h$ on $M$.\hfill\break
$\sec(p)(v,w)$ is the sectional curvature of the plane spanned by the
linearly independent vectors $v,w$ at $p$.\hfill\break
$\sec \geq k$ means that the sectional curvature of every plane at every point
is bounded from below by $k$.\hfill\break
$\curlR$ denotes the curvature operator.\hfill\break
$\curlR \geq c$ means that the eigenvalues of the curvature operator are bigger than
or equal to $c$ at every point on the manifold.\hfill\break
$\Gamma(h)^k_{ij}$ or $\uph\Gamma^k_{ij}$ refer to the Christoffel symbols 
of the metric $h$ in the coordinates $\{ x^k\},$
$$ \uph\Gamma^k_{ij} = {1 \on 2} h^{kl}( \parti{h_{il}}{x^j} + \parti{h_{jl}}{x^i}-   \parti{h_{ij}}{x^l} ).$$
For a diffeomorphism $F:M \to N$ we will sometimes consider $dF,$
a 1-form along $F$, defined by
$$dF(x) := \parti{F^{\al}}{x^k}dx^k({x}) \parti{}{y^\al}|_{(F(x))}.$$
For a general 1-form $\omega$ along $F,$ $\omega = \omega^{\al}_{i}(x)dx^i(x)\otimes  \parti{}{y^\al}|_{(F(x))},$ 
we define the norm of $\omega$
 with respect to $l$ (a metric on $M$) and $\ga$ (a metric on $N$) by
$$^{l,\ga}|\omega|^2(x) =  l^{ij}(x)\ga_{\al \be}(F(x)) \omega^{\al}_{i}(x) \omega^{\be}_{j}(x).$$
For example,
$${}^{l,\ga}|dF|^2(x) =   l^{ij}(x)\ga_{\al \be}(F(x)) \parti{F^{\al}}{x^i}(x) \parti{F^{\be}}{x^j} (x).$$  
\hfill\break
We define  $^{g,h}\grad dF,$ a ${0 \choose 2}$ tensor along $F$, by
$$ (^{g,h}\grad dF)_{ij}^{\al} := 
\Big( 
 \frac{\partial^2 F^{\al} }{\partial x^i \partial x^j}  
- \Gamma^{k}_{i j}(g) \parti{F^{\al}}{x^{k}} 
+ \Gamma^{\al}_{\be \si}(h)  \parti{F^{\be}}{x^{i}} 
\parti{F^{\si}}{x^{j}} \Big).$$\hfill\break
For a general  ${0 \choose 2}$ tensor $\psi$ along $F,$ 
$\psi = \psi^{\al}_{ij}(x) dx^i(x)\otimes  dx^j(x) \otimes  \parti{}{y^\al}|_{(F(x))},$ 
we define the norm of $\psi$
 with respect to $l$ (a metric on $M$) and $\ga$ (a metric on $N$) by
$$^{l,\ga}|\psi|^2 = \ga_{\al \be}(F(x)) l^{ks}(x) l^{ij}(x)\eta^{\al}_{ik}(x)\eta^{\be}_{js}(x).$$
For example
\begin{eqnarray*}
^{l,\ga}|^{g,h}\grad dF|^2  
&=& \ga_{\al \be}(F(x)) l^{ks}(x) l^{ij}(x) \Big( 
 \frac{\partial^2 F^{\al} }{\partial x^i \partial x^k}  
- \Gamma^{r}_{i k}(g) \parti{F^{\al}}{x^{r}} 
+ \Gamma^{\al}_{\eta \si}(h)  \parti{F^{\eta}}{x^{i}} 
\parti{F^{\si}}{x^{k}} \Big) \\ 
&&\Big( 
 \frac{\partial^2 F^{\be} }{\partial x^j \partial x^s}  
- \Gamma^{r}_{j s}(g) \parti{F^{\be}}{x^{r}} 
+ \Gamma^{\be}_{\phi \rho}(h)  \parti{F^{\phi}}{x^{j}} 
\parti{F^{\rho}}{x^{s}} \Big).
\end{eqnarray*}

\vskip 0.1 true in
\section*{Acknowledgements\markboth{Acknowledgements}{Acknowledgements}}

We would like to thank Peter Topping for helpful discussions
on Harmonic map heat flow and the Pseudolocality result of Perelman. 
Thanks to Klaus Ecker, Gerhard Huisken and Ernst Kuwert for their interest in and support of this
work.

\vskip 0.1 true in
\rm

\vskip 0.3 true in
\centerline{Mathematisches Institut, Eckerstr. 1, 79104 Freiburg im Br., Germany }
 \centerline{e-mail: msimon@mathematik.uni-freiburg.de}




\begin{thebibliography}{10}
\bibitem[1]{AnKn1}
Angenent, Sigurd; Knopf, Dan 
\newblock
An example of neckpinching for Ricci flow on $S\sp {n+1}$. 
\newblock Math. Res. Lett.  11  (2004),  no. 4, 493--518.

\bibitem[2]{AnKn2}
Angenent, Sigurd; Knopf, Dan 
\newblock
Precise asymptotics of the Ricci flow neckpinch . 
\newblock Pre-print

\bibitem[3]{Be}
Berger, M.
\newblock Les Varietes Riemanniennes a courbure positive, 
\newblock Bull.Soc.Math.France 87, (1960).
\bibitem[4]{CCCY} 
Cao,H.D., Chow, B., Chu,S.C., Yau,S.T.
\newblock Collected papers on the Ricci flow
\newblock Series in Geometry and Topology, Vol.37, International Press.

\bibitem[5]{BuBuIv}
Burago, D., Burago, Y., Ivanov, S.
\newblock A course in Metric Geometry, 
\newblock Graduate studies in Math., Vol. 33, American Math. Soc.


\bibitem[6]{BGP}
Burago, Yu., Gromov, M., Perelman,G.
\newblock A.D. Alexandrov spaces with curvature bounded below.
\newblock  Russian Math surveys, Vol. 47, pp. 1-58 (1992).
\bibitem[7]{ChKn}
Chow,B. Knopf,D.
\newblock The Ricci Flow: An Introduction
\newblock Math. Surveys and Mono. Volume 110, American Math. Soc.


\bibitem[8]{ChCo96}
Cheeger,J., Colding,T,.
\newblock Lower bounds on the Ricci curvature and the almost rigidity of warped
products
\newblock Ann. of Math., 144 , pp. 189-237, (1996)





\bibitem[8]{ChCo}
Cheeger,J., Colding,T,.
\newblock On the structure of spaces with Ricci with curvature bounded from below I,
\newblock J.differential Geometry 46, 406-480 (1997)


\bibitem[9]{CGT}
Cheeger,J., Gromov,M., Taylor,M.
\newblock Finite propagation speed, kernel estimates for functions of the Laplace operator,
 and the geometry of complete Riemannian manifolds,
\newblock {\em J.Differential Geom.},  17,no. 1, 15 -- 53, (1982).

\bibitem[10]{ChEb}
Cheeger, J., Ebin,D.
\newblock Comparison Theorems in Riemannian Geometry
\newblock {\em North Holland Publishing Company}, 1975

\bibitem[11]{De}
DeTurck,~D.
\newblock  Deforming metrics in the direction of their
 Ricci tensors,
\newblock {\em J.Differential Geom.},  18,no. 1, 157 -- 162, (1983).

\bibitem[12]{FuYa}
Fukaya, K., Yamaguchi, T.
\newblock  The fundamental groups of almost non-negatively curved manifolds
\newblock {\em Ann. of Math. }, (2), 136, 253-333, (1992).


\bibitem[13]{Ha82}
Hamilton,~R.S.
\newblock Three manifolds with positive 
Ricci-curvature
\newblock {\em J.Differential Geom.},  17, no. 2, 255 -- 307, (1982).

\bibitem[14]{Ha86}
Hamilton,~R.S.
\newblock Four manifolds with positive curvature operator
 \newblock {\em J. Differential
Geom} 24 no. 2 , 153 -- 179, (1986).

\bibitem[15]{Ha93}
Hamilton,~R.S.
\newblock
 Eternal solutions to the Ricci flow,
\newblock {\em Journal of Diff. Geom.}, 38, 1-11 (1993)

\bibitem[16]{Ha95}
Hamilton,~R.S.
\newblock
 The formation of singularities in the Ricci flow,
\newblock {\em Collection: Surveys in differential geometry}, 
Vol. II (Cambridge, MA), 7--136, (1995).

\bibitem[17]{Ha95a}
Hamilton,~R.S.
\newblock A compactness property of the Ricci Flow
\newblock {\em American Journal of Mathematics}, 117, 545--572, (1995)


\bibitem[18]{Ha99}
Hamilton,~R.S.
\newblock Non-Singular solutions of the Ricci Flow on Three-Manifolds
\newblock {\em Comm. Anal. Geom.}, vol 7., no. 4, 695--729, (1999)

\bibitem[19]{He}
Hebey, E.
\newblock Nonlinear Analysis on Manifolds: Sobolev Spaces and Inequalities
\newblock Courant Lecture Notes in Mathematics, New York University, 1999


\bibitem[20]{Hi}
Hirsch, M.W.
\newblock Differential Topology
\newblock {\em Springer,1976}.

\bibitem[21]{Iv}
Ivey, T. 
\newblock Ricci solitons on compact three manifolds, 
\newblock Diff.Geom.Appl. 3, 301--307 (1993)


\bibitem[22]{Kl}
Klingenberg, W.
\newblock Contributions to Riemannian geometry in the large
\newblock  Ann.Math. 69, (1959)

\bibitem[23]{Lo} 
Lohkamp, Jo.
\newblock Negatively Ricci curved manifolds.  
\newblock Bull. Amer. Math. Soc. (N.S.)  27  (1992),  no. 2, 288--291

\bibitem[24]{LSU}
Ladyzenskaja,O.A., Solonnikov, V.A., Uralceva, N.N.
 \newblock Linear and quasilinear equations of parabolic type,
\newblock {\em
Transl.Amer.Math.Soc.,1986 }

\bibitem[25]{Pe1}
Perelman,G.,
\newblock 
The entropy formula for the Ricci flow and its geometric applications
\newblock MarthArxiv link: math.DG/0211159

\bibitem[26]{Pe2}
Perelman,G.,
\newblock 
Ricci flow with surgery on three manifolds
\newblock MarthArxiv link: math.math.DG/0303109


\newblock {\em Springer,1991}.
\bibitem[27]{Pet}
Peters, S.
\newblock Konvergenz Riemannscher Mannigfaltigkeiten.
\newblock {\em Bonner Math. Schr.} 169, (1986)


\bibitem[28]{Pe}
Peterson,P.,
\newblock Riemannian Geometry
\newblock {\em Springer,1991}.

\bibitem[29]{Ra}
Rauch, H.E.
\newblock A contribution to differential geometry in the large
\newblock Ann. Math. 54, 38--55 (1951)


\bibitem[30]{Si} 
Simon, M.,
\newblock Deformation of $C^0$ Riemannian metrics in the
direction of their Ricci curvature,
\newblock {\em Comm. Anal. Geom.},
10, no. 5, 1033-1074, (2002)

\bibitem[31]{Si2} 
Simon, M.,
\newblock
A class of Riemannian manifolds that pinch when evolved by Ricci flow.  
\newblock Manuscripta Math.  101  (2000),  no. 1, 89--114.

\bibitem[32]{Ya} 
Yamaguchi, T.,
\newblock Collapsing and pinching under a lower curvature bound
\newblock {\em Annals of mathematics},
133, 317-357 (1999)

\bibitem[33]{ShYa00} 
Shioya, T., Yamaguchi, T.,
\newblock Collapsing three-manifolds under a lower curvature bound.  
\newblock J. Differential Geom.  56,  no. 1, 1--66.  (2000)


\bibitem[34]{ShYa05} 
Shiyoa, T., Yamaguchi, T.,
\newblock  Volume collapsed three-manifolds with a lower curvature bound. 
\newblock  Math. Ann.  333,  no. 1, 131--155. (2005)


\bibitem[35]{Sh}
Shi, Wan-Xiong.,
\newblock Deforming the metric
on complete Riemannian manifolds 
\newblock {\em J.Differential Geometry} , 30, 223--301,(1989).
\end{thebibliography}
\end{document}